\documentclass[12pt]{article}

\newcommand{\slaover}[1]{}

\usepackage{standalone}

\usepackage{fullpage}
\usepackage{amssymb, amsmath, amsthm}
\usepackage{color}
\usepackage{standalone}
\usepackage{hyperref}
\renewcommand{\d}{\,\mathrm{d}}

\newcommand{\dx}{\,\mathrm{d}x}
\newcommand{\dy}{\,\mathrm{d}y}

\newcommand{\connectionrule}{$(x,y) \sim (x',y') \iff |x-x'|_{\pi e^{R/2}} \leq e^{\frac12(y+y')}$}

\newtheorem{thm}{Theorem}
\newtheorem{lem}[thm]{Lemma}
\newtheorem{cor}[thm]{Corollary}
\newtheorem{proposition}[thm]{Proposition}

\newcommand{\blokje}{\hfill $\blacksquare$\\}
\renewcommand{\theta}{\vartheta}
\renewcommand{\phi}{\varphi}
\renewcommand{\epsilon}{\varepsilon}
\renewcommand{\rho}{\varrho}
\usepackage{enumerate}
\usepackage[font=small]{caption}
\usepackage{pgf,tikz}
\usepackage{ifthen}
\usetikzlibrary{arrows,patterns}
\usepackage{etoolbox}

\newcommand{\Acal}[0]{\ensuremath{{\mathcal A}}}
\newcommand{\Bcal}[0]{\ensuremath{{\mathcal B}}}

\newcommand{\Ecal}[0]{\ensuremath{{\mathcal E}}}

\newcommand{\Gtil}[0]{\tilde{G}}

\newcommand{\Vtil}[0]{\tilde{V}}
\newcommand{\Wtil}[0]{\tilde{W}}

\newcommand{\elltil}[0]{\tilde{\ell}}

\newcommand{\Pee}[0]{\ensuremath{{\mathbb P}}}

 \newcommand{\eps}{\varepsilon}

\DeclareMathOperator{\Po}{Po}

\DeclareMathOperator{\dd}{d}

\newcommand{\Gammatil}[0]{\ensuremath{\tilde{\Gamma}}}

\renewenvironment{proof}{\vspace{1ex}\noindent{\bf Proof:}}{\hspace*{\fill}$\blacksquare$\vspace{1ex}}
\newenvironment{proofof}[1]{\vspace{1ex}\noindent{\bf Proof of #1:}}{\hspace*{\fill}$\blacksquare$\vspace{1ex}}

\begin{document}

\title{The diameter of KPKVB random graphs}
\author{Tobias M\"uller\thanks{Bernoulli Institute, Groningen University, {\tt tobias.muller@rug.nl}. 
Supported in part by the Netherlands Organisation for Scientific Research (NWO) under project nos 612.001.409 and 639.032.529} 
\and Merlijn Staps\thanks{Department of Ecology and Evolutionary Biology, Princeton University, {\tt merlijnstaps@gmail.com}. This paper is the result of this 
author's MSc thesis research project, carried out at Utrecht Uniersity under the supervision of the first author. 
The thesis is available from~\cite{MerlijnThesis}.}}
\date{\today}
\maketitle

\begin{abstract}
We consider a random graph model that was recently proposed as a model for complex networks 
by Krioukov et al.~\cite{Krioukov}.
In this model, nodes are chosen randomly inside a disk in the hyperbolic plane and two nodes are connected if they are at most a certain 
hyperbolic distance from each other. It has been previously shown that this model has various properties associated with complex networks, including 
a power-law degree distribution and a strictly positive clustering coefficient. 
The model is specified using three parameters: the number of nodes $N$, which we think of as going to infinity, and $\alpha, \nu > 0$, which we think of 
as constant. Roughly speaking $\alpha$ controls the power law exponent of the degree sequence and $\nu$ the average degree.

Earlier work of Kiwi and Mitsche~\cite{Kiwi} has shown that when $\alpha < 1$ (which corresponds to the exponent of the power law degree sequence 
being $< 3$) then the diameter of the largest component is a.a.s.~at most polylogarithmic 
in $N$. Friedrich and Krohmer~\cite{FK} have shown it is a.a.s.~$\Omega(\log N)$ and they improved the exponent of the 
polynomial in $\log N$ in the upper bound.
Here we show the maximum diameter over all components is a.a.s.~$O(\log N)$ thus giving a bound that is tight up to a multiplicative constant.
\end{abstract}

\section{Introduction}

The term \emph{complex networks} usually refers to various large real-world networks, occurring diverse fields of science, that
appear to exhibit very similar graph theoretical properties. 
These include having a constant average degree, a so-called power-law degree sequence, clustering and ``small distances''.
In this paper we will study a random graph model that was recently proposed as a model for complex networks and has the above properties.
We refer to it as the Krioukov-Papadopoulos-Kitsak-Vahdat-Bogu\~n\'a model, or KPKVB model, after its inventors~\cite{Krioukov}. 
We should however maybe point out that many authors simply refer to the model as ``hyperbolic random geometric graphs'' or even 
``hyperbolic random graphs''.
In the KPKVB model a random geometric graph is constructed in the hyperbolic plane. We use the Poincar\'e disk representation of the hyperbolic 
plane, which is obtained when the unit disk 
${\mathbb D} = \{(x,y) \in \mathbb{R}^2 : x^2 + y^2 < 1\}$ is equipped with the metric given by the differential form 
$\d s^2 = 4 \frac{\d x^2 + \d y^2}{(1-x^2-y^2)^2}$. (This means that the length of a curve $\gamma : [0,1]\to{\mathbb D}$ under the metric is given by 
$2 \int_0^1\frac{\sqrt{(\gamma_1'(t))^2+(\gamma_2'(t))^2}}{1-\gamma_1^2(t)-\gamma_2^2(t)}{\dd}t$.) 
For an extensive, readable introduction to hyperbolic geometry and the various models and properties of the hyperbolic plane, the reader
could consult the book of Stillwell~\cite{Stillwell}.
Throughout this paper we will represent points in the hyperbolic plane by polar coordinates $(r,\theta)$, 
where $r \in [0,\infty)$ denotes the hyperbolic distance of a point to the origin, and $\theta$ denotes its angle with the positive $x$-axis. 

We now discuss the construction of the KPKVB random graph. The model has three parameters: the number of vertices $N$ 
and two additional parameters $\alpha, \nu > 0$. 
Usually the behavior of the random graph is studied for $N \to \infty$ for a fixed choice of $\alpha$ and $\nu$. 
We start by setting $R = 2\log(N/\nu)$. Inside the disk $\mathcal{D}_R$ of radius $R$ centered at the origin in the hyperbolic plane 
we select $N$ points, independent from each other, according to the probability density $f$ on $[0,R] \times (-\pi, \pi]$ given by
\[
f(r,\theta) = \frac{1}{2\pi} \frac{\alpha \sinh(\alpha r)}{\cosh(\alpha R)-1} \qquad ( r \in [0,R], \theta \in (-\pi,\pi]).
\]
We call this distribution the $(\alpha,R)$-quasi uniform distribution. For $\alpha=1$ this corresponds to the uniform distribution on 
$\mathcal{D}_R$. 
We connect points if and only if their hyperbolic distance is at most $R$. In other words, two points are connected if their hyperbolic distance 
is at most the (hyperbolic) radius of the disk that the graph lives on. We denote the random graph we have thus obtained by 
$G(N;\alpha,\nu)$.

As observed by Krioukov et al.~\cite{Krioukov} and rigorously shown by Gugelmann et al.~\cite{Gugelmann}, the degree distribution
follows a power law with exponent $2\alpha+1$, the average degree tends to $2\alpha^2 \nu / \pi(\alpha-1/2)^2$ when $\alpha > 1/2$, 
and the (local) clustering coefficient is bounded away from zero a.a.s.
(Here and in the rest of the paper a.a.s.~stands for {\em asymptotically almost surely}, meaning with probability tending to one as $N\to\infty$.)
Earlier work of the first author with Bode and Fountoulakis~\cite{BFM_giant} and with Fountoulakis~\cite{FM} has established the 
``threshold for a giant component'': when $\alpha < 1$ then there always is a unique component of size linear in $N$ no matter how small 
$\nu$ (and hence the average degree) is; when $\alpha > 1$ all components are sublinear no matter the value of $\nu$; and when $\alpha=1$ then 
there is a critical value $\nu_{\text{c}}$ such that for $\nu < \nu_{\text{c}}$ all components are sublinear and for $\nu > \nu_{\text{c}}$ there 
is a unique linear sized component (all of these statements holding a.a.s.). 
Whether or not there is a giant component when $\alpha=1$ and $\nu=\nu_{\text{c}}$ remains an open problem.

In another paper of the first author with Bode and Fountoulakis~\cite{BFM_connected} it was shown that $\alpha=1/2$ is the threshold for connectivity:
for $\alpha < 1/2$ the graph is a.a.s.~connected, for $\alpha>1/2$ the graph is a.a.s.~disconnected, and when $\alpha=1/2$ the probability 
of being connected tends to a continuous, nondecreasing function of $\nu$ which is identically one for $\nu \geq \pi$ and strictly less than one 
for $\nu < \pi$.

Friedrich and Krohmer~\cite{FK_clique} studied the size of the largest clique as well as the number of cliques of a given size.
Bogu\~{n}a et al.~\cite{sustaining} and Bl\"asius et al.~\cite{BlasiusMLE} considered fitting the KPKVB model to data using maximum likelihood 
estimation.
Kiwi and Mitsche~\cite{KiwiSpectral} studied the spectral gap and related properties, and Bl\"asius et al.~\cite{BlasiusTW} considered the treewidth 
and related parameters of the KPKVB model.

Abdullah et al.~\cite{ABF} considered typical distances in the graph. That is, they sampled two vertices of the graph uniformly 
at random from the set of all vertices and consider the (graph-theoretic) distance between them. They showed that this distance between two 
random vertices, conditional on the two points falling in the same component, is precisely $(c + o(1)) \cdot \log\log N$ 
a.a.s. for $1/2 < \alpha < 1$, where $c := - 2 \log(2\alpha-1)$.

Here we will study another natural notion related to the distances in the graph, the graph diameter.
Recall that the \emph{diameter} of a graph $G$ is the supremum of the graph distance $d_G(u,v)$ over all pairs $u$, $v$ of vertices (so it is infinite if 
the graph is disconnected).
It has been shown previously by Kiwi and Mitsche~\cite{Kiwi} that for $\alpha \in (\frac12,1)$ the largest component of $G(N;\alpha,\nu)$ has a diameter that is 
$O\left( (\log N)^{8/(1-\alpha)(2-\alpha)} \right)$ a.a.s.
This was subsequently improved by Friedrich and Krohmer~\cite{FK} to $O\left( (\log N)^{1/(1-\alpha)} \right)$.
Friedrich and Krohmer~\cite{FK} also gave an a.a.s.~lower bound of $\Omega(\log N)$.
We point out that in these upper bounds the exponent of $\log N$ tends to infinity as $\alpha$ approaches one.

Here we are able to improve the upper bound to $O(\log N)$, which is sharp up to a multiplicative constant.
We are able to prove this upper bound not only in the case when $\alpha < 1$ but also in the case when $\alpha=1$ and $\nu$ is sufficiently large.

\begin{thm}\label{thm:main}
Let $\alpha, \nu > 0$ be fixed. If either 
\begin{enumerate}[(i)]
\item $\frac12 < \alpha < 1$ and $\nu > 0$ is arbitrary, or; 
\item $\alpha = 1$ and $\nu$ is sufficiently large, 
\end{enumerate}
then, a.a.s.~as $N\to\infty$, every component of $G(N; \alpha, \nu)$ has diameter $O(\log(N))$.
\end{thm}

We remark that our result still leaves open what happens for other choices of $\alpha, \nu$ as well as several related questions.
See Section~\ref{sec:conclusion} for a more elaborate discussion of these.

\subsection{Organization of the paper}

In our proofs we will also consider a Poissonized version of the KPKVB model, where the number of points is not fixed but is sampled from a Poisson 
distribution with mean $N$. This model is denoted $G_{\Po}(N; \alpha, \nu)$. It is convenient to work with this Poissonized version 
of the model as it has the advantage that the numbers of points in disjoint regions are independent (see for instance~\cite{Kingman}).

The paper is organized as follows. In Section \ref{sec:idealized} we discuss a somewhat simpler random geometric graph $\Gamma$, introduced 
in~\cite{FM}, that behaves approximately 
the same as the (Poissonized) KPKVB model. The graph $\Gamma$ is embedded into a rectangular 
domain $\mathcal{E}_R$ in the Euclidean plane $\mathbb{R}^2$. In Section \ref{sec:boxes} we discretize this simplified model by dissecting $\mathcal{E}_R$ into 
small rectangles. In Section \ref{sec:paths} we show how to construct short paths in $\Gamma$. 
The constructed paths have length $O(\log(N))$ unless there exist large regions that do not contain any vertex of $\Gamma$. 
In Section \ref{sec:bounding} we use the observations of Section \ref{sec:paths} to formulate sufficient conditions for the components of the graph $\Gamma$ to have diameter $O(\log N)$. In Section \ref{sec:probbounds} we show that the probability that $\Gamma$ fails to satisfy these conditions tends to $0$ as $N \to \infty$. We also translate these results to the KPKVB model, and combine everything into a proof of Theorem \ref{thm:main}.

\section{The idealized model\label{sec:idealized}}

We start by introducing a somewhat simpler random geometric graph, introduced in~\cite{FM}, that will be used as an approximation of the KPKVB model.
Let $X_1$, $X_2$, \ldots $\in \mathcal{D}_R$ be an infinite supply of points chosen according the $(\alpha,R)$-quasi uniform distribution on 
$\mathcal{D}_R$ described above. Let $G = G(N; \alpha,\nu)$ and $G_{\Po} = G_{\Po}(N; \alpha, \nu)$. Let $Z \sim \Po(N)$ be the number of 
vertices of $G_{\Po}$. By taking $\{X_1,\ldots,X_N\}$ as the vertex set of $G$ and $\{X_1,\ldots,X_Z\}$ as the vertex set of $G_{\Po}$, we 
obtain a coupling between $G$ and $G_{\Po}$.

\begin{figure}[h!]
\centering
\begin{tikzpicture}[scale=1.3]

\def\xmin{2.6}
\def\xmax{5.6}
\def\ymin{-1}
\def\ymax{1}
\pgfmathsetmacro{\halfwidth}{(\xmax-\xmin)/2}

\draw (0,0) circle (1);
\node at (0,1.3) {$\mathcal{D}_R$};

\draw [->] (1.2,0.05) -- (2.4,0.05);
\node at (1.8,0.35) {$\Psi$};

\draw (\xmin,\ymin) -- (\xmin,\ymax) -- (\xmax,\ymax) -- (\xmax,\ymin) -- (\xmin,\ymin);
\node at (\xmin+\halfwidth,\ymax+0.3) {$\mathcal{E}_R$};
\node at (\xmin,\ymin-0.25) {\scriptsize $-\frac{\pi}{2} e^{R/2}$};
\node at (\xmax,\ymin-0.25) {\scriptsize $\frac{\pi}{2} e^{R/2}$};

\draw [<->] (\xmax+0.2,\ymin) -- (\xmax+0.2,\ymax);
\node at (\xmax+0.35,0) {\scriptsize $R$};

\draw [gray,dashed] (-1,0) -- (1,0);
\draw [gray,dashed] (0,-1) -- (0,1);
\draw [<->] (-0.03,-0.03) -- (-0.57,-0.78);
\node at (-0.45,-0.35) {\scriptsize $R$};

\draw [gray,dashed] (\xmin+\halfwidth,\ymin) -- (\xmin+\halfwidth,\ymax);
\draw [gray,dashed] (\xmin+0.5*\halfwidth,\ymin) -- (\xmin+0.5*\halfwidth,\ymax);
\draw [gray,dashed] (\xmin+1.5*\halfwidth,\ymin) -- (\xmin+1.5*\halfwidth,\ymax);

\draw [fill=red,red] (0.2,-0.3) circle [radius=0.05];
\draw [fill=red,red] (\xmin + \halfwidth -0.98*\halfwidth/3.141593,0.28) circle [radius=0.05];

\draw [fill=blue,blue] (-0.1,0.1) circle [radius=0.05];
\draw [fill=blue,blue] (\xmin + \halfwidth + 3*\halfwidth/4,0.72) circle [radius=0.05];

\draw [fill=green,green] (0.5,0.7) circle [radius=0.05];
\draw [fill=green,green] (\xmin + \halfwidth + 0.9505468*\halfwidth/3.141593,-0.72) circle [radius=0.05];

\end{tikzpicture}
\caption{$\Psi$ maps $\mathcal{D}_R$ to a rectangle $\mathcal{E}_R \subset \mathbb{R}^2$.\label{fig:psi}}
\end{figure}
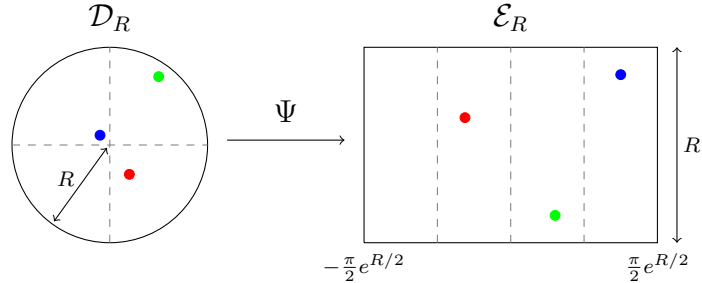

We will compare our hyperbolic random graph to a random geometric graph that lives on the Euclidean plane. To this end, we introduce the 
map $\Psi: \mathcal{D}_R \to \mathbb{R}^2$ given by
$
\Psi: (r,\theta) \mapsto \left(\theta \cdot \frac12 e^{R/2}, R-r\right). 
$
The map $\Psi$ works by taking the distance of a point to the boundary of the disk as $y$-coordinate and the angle of the point as $x$-coordinate 
(after scaling by $\frac12 e^{R/2}$) 
The image of $\mathcal{D}_R$ under $\Psi$ is the rectangle 
$\mathcal{E}_R = (-\frac{\pi}{2} e^{R/2}, \frac{\pi}{2} e^{R/2}] \times [0,R] \subset \mathbb{R}^2$ (Figure \ref{fig:psi}). 

On $\mathcal{E}_R$ we consider the Poisson point process $\mathcal{P}_{\alpha, \lambda}$ with intensity function $f_{\alpha,\lambda}$ defined 
by $f_{\alpha,\lambda}(x,y) = \lambda e^{-\alpha y}$. 
We will denote by $V_{\alpha,\lambda}$ the point set of this Poisson process. 
We also introduce the graph $\Gamma_{\alpha,\lambda}$, with vertex set $V_{\alpha,\lambda}$, where points $(x,y), (x',y') \in V_{\alpha,\lambda}$ 
are connected if and only if $|x-x'|_{\pi e^{R/2}} \leq e^{\frac12(y+y')}$.  
Here $|a-b|_{d} = \inf_{k \in \mathbb{Z}} |a-b+kd|$ denotes the distance between $a$ and $b$ 
modulo $d$.

If we choose $\lambda = \frac{\nu \alpha}{\pi}$ it turns out that $V_{\lambda}$ can be coupled to the image of the vertex set of $G_{\Po}$ 
under $\Psi$ and that the connection rule of $\Gamma_{\lambda}$ approximates the connection rule of $G_{\Po}$. In particular, we have the following:

\begin{lem}[\cite{FM}, Lemma 27]\label{lem:coupling}
Let $\alpha>\frac12$. There exists a coupling such that a.a.s.  $V_{\alpha, \nu \alpha/\pi}$ is the image of the vertex set of $G_{\Po}$ under $\Psi$. 
\end{lem}

\noindent
Let $\tilde{X}_1$, $\tilde{X}_2$, $\ldots \in \mathcal{E}_R$ be the images of $X_1$, $X_2$, \ldots\ under $\Psi$. On the coupling space of 
Lemma~\ref{lem:coupling}, a.a.s.\ we have $V_{\lambda} = \{\tilde{X}_1, \ldots, \tilde{X}_Z\}$.

\begin{lem}[\cite{FM}, Lemma 30]\label{lem:coupling-edges}
Let $\alpha>\frac12$. On the coupling space of Lemma \ref{lem:coupling}, a.a.s.\ it holds for $1 \leq i,j \leq Z$ that\index{coupling}
\begin{itemize}
\item[(i)] if $r_i,r_j \geq \frac12R$ and $\tilde{X}_i\tilde{X}_j \in E(\Gamma_{\alpha,\nu \alpha/\pi})$, then $X_iX_j \in E(G_{\Po})$.
\item[(ii)] if $r_i,r_j \geq \frac34R$, then $\tilde{X}_i\tilde{X}_j \in E(\Gamma_{\alpha,\nu \alpha/\pi}) \iff X_iX_j \in E(G_{\Po})$. 
\end{itemize} Here $r_i$ and $r_j$ denote the radial coordinates 
of $X_i,X_j \in \mathcal{D}_R$.\index{coupling!between KPKVB model and idealized model}
\end{lem}

Lemma \ref{lem:coupling-edges} will prove useful later because as it turns out cases (i) and (ii) cover almost all the edges in the graph.

For $(A_i)_i, (B_i)_i$ two sequences of events with $A_i$ and $B_i$ defined on the same probability 
space $(\Omega_i,\Acal_i,\Pee_i)$, we 
say that $A_i$ happens a.a.s.\ conditional on $B_i$ if $\mathbb{P}(A_i \mid B_i) \to 1$ as $i \to \infty$.
By a straightforward adaptation of the proofs given in~\cite{FM}, it can be shown that also:

\begin{cor}\label{cor:koppel}
The conclusions of Lemmas~\ref{lem:coupling} and~\ref{lem:coupling-edges} also hold conditional on the event $Z=N$. 
\end{cor}

\begin{figure}[h!]\centering
\includegraphics[height=0.4\linewidth]{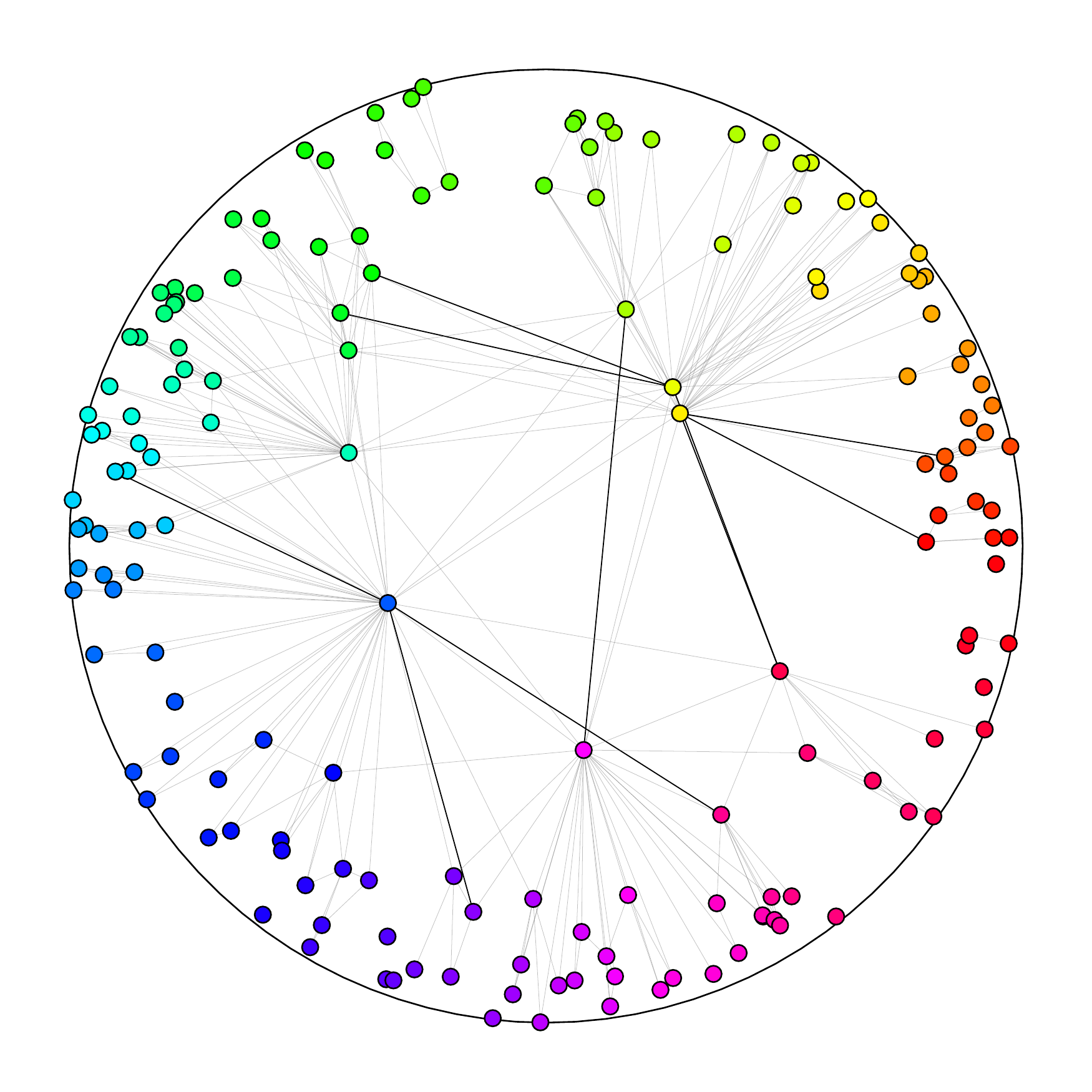} \includegraphics[height=0.4\linewidth]{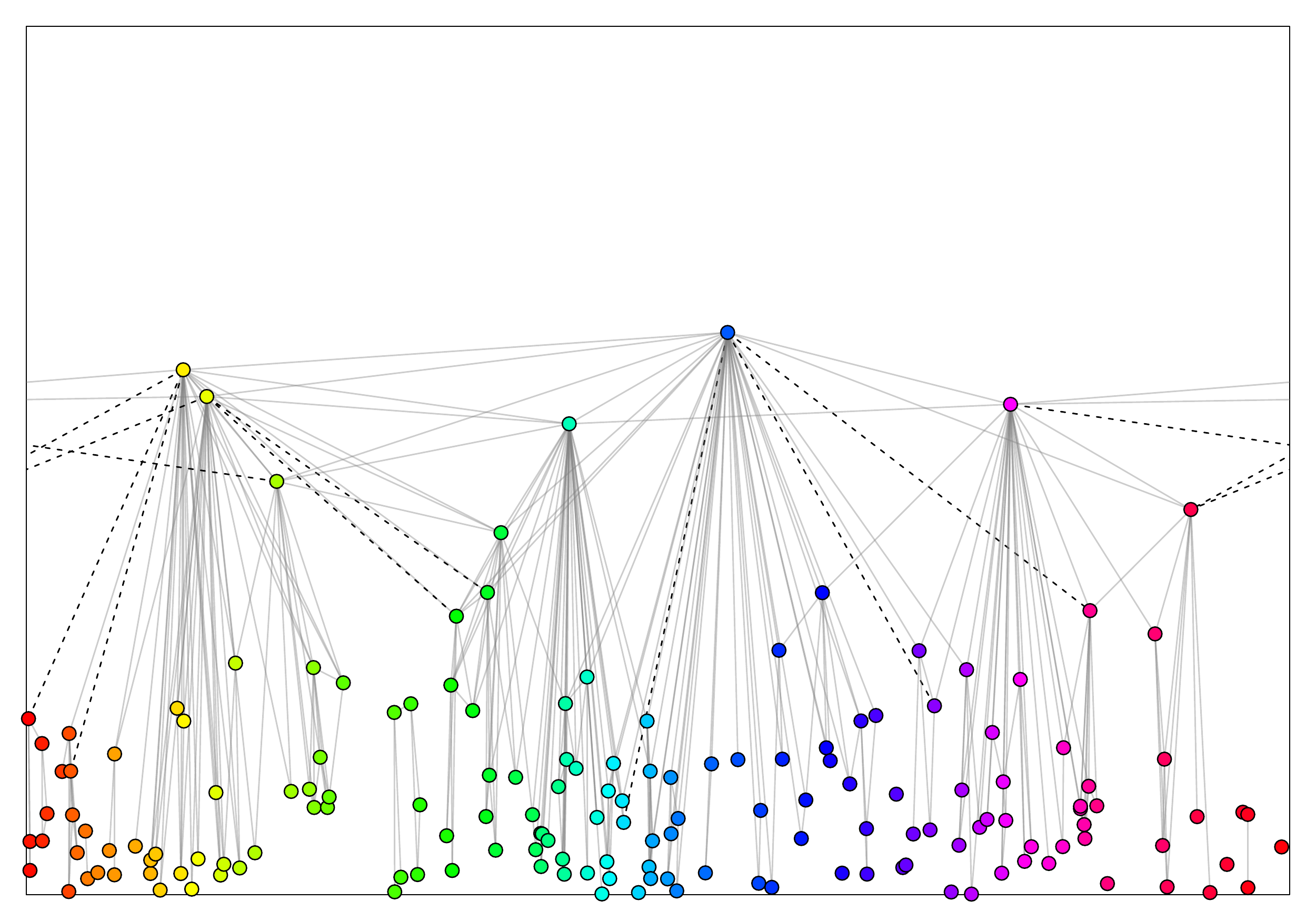}
\caption{An example of the Poissonized KPKVB random graph $G_{\Po}$ (left) and the graph $\Gamma_{\alpha,\nu \alpha/\pi}$ (right), under 
the coupling of Lemma \ref{lem:coupling}. 
The graph $G_{\Po}$ is drawn in the native model of the hyperbolic plane\index{hyperbolic plane!native representation}, where a point with 
hyperbolic polar coordinates $(r,\theta)$ is plotted with Euclidean polar coordinates $(r,\theta)$. Points are colored based on their 
angular coordinate. 
The edges for which the coupling fails are drawn in black 
in the picture of $G_{\Po}$ and as dotted lines in the picture of $\Gamma_{\nu \alpha/\pi}$. 
The parameters used are $N=200$, $\alpha=0.8$ and $\nu = 1.3$. \label{fig:DRER}}
\end{figure}

\noindent
In other words, the corollary states that the probability that the conclusions of Lemmas~\ref{lem:coupling} and~\ref{lem:coupling-edges}
fail, given that $Z=N$, is also $o(1)$.
For completeness, we prove this as Lemmas~\ref{lem:coupling-2} and~\ref{lem:coupling-edges-2} in the appendix. 
An example of $G_{\Po}$ and $\Gamma_{\nu \alpha/\pi}$ is shown in Figure \ref{fig:DRER}.

\section{Deterministic bounds}

For the moment, we continue in a somewhat more general setting, where  $V \subset \mathcal{E}_R$ is any finite set of points and $\Gamma$ is the graph 
with vertex set $V$ and connection rule \connectionrule.

\subsection{A discretization of the model~\label{sec:boxes}}

 We dissect $\mathcal{E}_R$ into a number of rectangles, which have the property 
that vertices of $\Gamma$ in rectangles with nonempty intersection are necessarily connected by an edge. This is done as follows. 
First, divide $\mathcal{E}_R$ into $\ell+1$ 
layers\index{layer} $L_0$, $L_1$, \ldots, $L_{\ell}$, where \index{$L_i$} \index{$\ell$ (number of layers)}
$$
L_i = \{(x,y) \in \mathcal{E}_R: i \log(2) \leq y < (i+1) \log(2)\}
$$ 
for $i < \ell$ and $L_{\ell} = \{(x,y) \in \mathcal{E}_R: y \geq \ell \log(2)\}$. Here $\ell$ is defined by 

\begin{equation}\label{eq:elldef} 
\ell := \left\lfloor \frac{\log(6\pi) + R/2}{\log(2)} \right\rfloor. 
\end{equation}

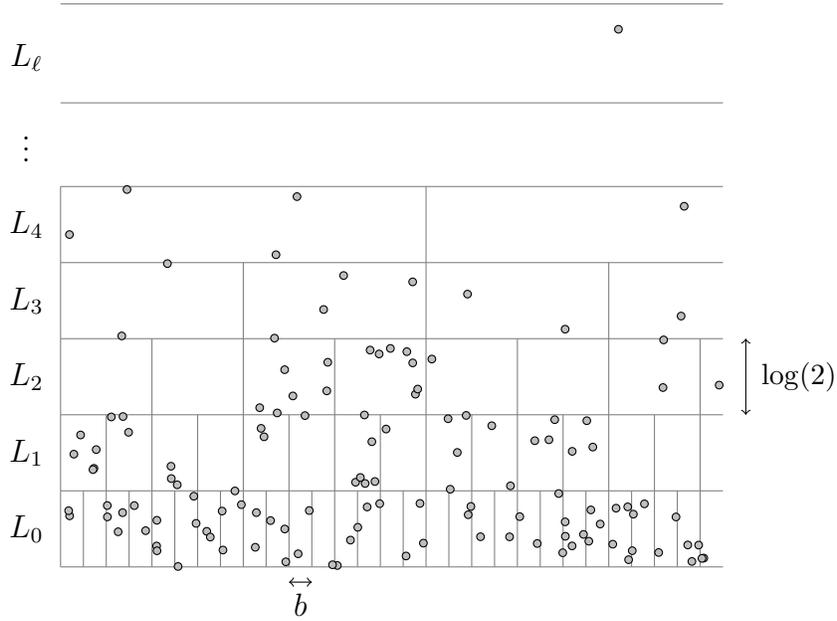
\begin{figure}[h!]
\centering
\resizebox{0.7\linewidth}{!}{
\begin{tikzpicture}
\def\xscale{0.3}
\def\columns{29}
\def\rows{5}
\def\n{8}
\pgfmathtruncatemacro{\rowsminusone}{\rows - 1}
\pgfmathtruncatemacro{\rowsplusone}{\rows + 1}
\pgfmathtruncatemacro{\rowsplustwo}{\rows + 2}
\pgfmathtruncatemacro{\columnsa}{\columns/2}


\foreach \x in {0,...,\rows}{
\draw [gray] (0, \x) -- (\columns*\xscale, \x);
}

\draw [gray] (0, \rows + 1.1) -- (\columns*\xscale, \rows+1.1);
\draw [gray] (0, \rows + 2.4) -- (\columns*\xscale, \rows+2.4);


\foreach \y in {0,...,\rowsminusone} {
    \pgfmathtruncatemacro{\columnsa}{\columns/2^\y}
	\foreach \x in {0,...,\columnsa}{
	\pgfmathsetmacro{\z}{2^\y*\x}
	\pgfmathtruncatemacro{\w}{\z}
	\ifthenelse{\w < \columns}{
	\draw [gray] (\z*\xscale,\y) -- (\z*\xscale,\y+1);
	}{}
	} 
}


{

\tikzset{%
	examplebox/.style = {black,fill=gray!50!white,line width=1.2pt},
    neighbor/.style={black,fill=gray!10!white,line width=0.8pt},
}

%

}


\foreach \y in {0, ..., \rowsminusone}{
\node at (-1.5*\xscale,\y+0.5) {$L_\y$};
}

\node at (-1.5*\xscale,\rows+0.6) {$\vdots$};
\node at (-1.5*\xscale,\rows+1.7) {$L_\ell$};

\def\bottomy{0}

\draw [<->] (10*\xscale,\bottomy-0.2) -- (11*\xscale,\bottomy-0.2);
\node at (10.5*\xscale,\bottomy-0.5) {$b$};

\draw [<->] (\columns*\xscale + 0.3, 2) -- (\columns*\xscale + 0.3, 3);
\node at (\columns*\xscale + 1,2.5) {\small $\log(2)$};


\foreach \r in {0,...,\rowsminusone}{
	\pgfmathtruncatemacro{\n}{2*2^\rows/2^\r}
	\foreach \t in {0,...,\n}{
	\pgfmathsetmacro{\randomx}{rnd*\columns}
	\pgfmathsetmacro{\randomy}{rnd}
	\draw [black,fill=gray!50!white] (\randomx*\xscale, \r+\randomy) circle [radius=0.05];
	}
}

	\pgfmathsetmacro{\randomx}{rnd*\columns}
	\pgfmathsetmacro{\randomy}{rnd}
	\draw [black,fill=gray!50!white] (\randomx*\xscale, \rowsplusone+0.2+\randomy) circle [radius=0.05];
	
%
%

\end{tikzpicture}
}
\caption{Partitioning $\mathcal{E}_R$ with boxes. All layers except $L_{\ell}$ have height $\log(2)$. 
The boxes in layer $L_i$ have width $2^i b$, where $b \in [\frac16,\frac13)$ is the width of a box in $L_0$. 
The small circles serve as an example of $V$.  
\label{fig:boxes}}
\end{figure}

\noindent
Note that this gives $6\pi e^{R/2} \geq 2^{\ell} > 3\pi e^{R/2}$. We divide $L_i$ into $2^{\ell-i}$ (closed) rectangles 
 of equal width $2^{i-\ell} \cdot \pi e^{R/2} = 2^{i} \cdot b$, where $b = 2^{-\ell} \cdot \pi e^{R/2} \in [\frac16,\frac13)$ is the 
 width of a rectangle in the lowest layer $L_0$ (Figure \ref{fig:boxes}). In each layer, one of the rectangles has its left edge on the line $x=0$. 
We have now partitioned $\mathcal{E}_R$ into $2^{\ell+1}-1 = O(N)$ boxes\index{box}. 

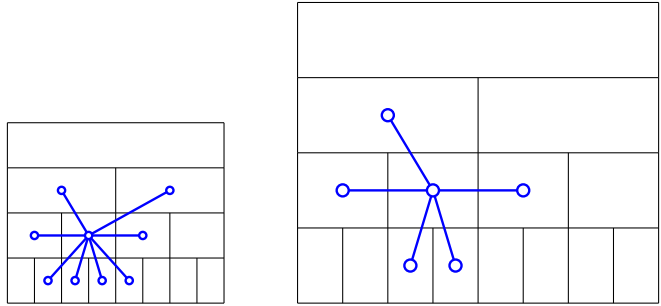
\begin{figure}[ht!] \centering
%
%
%
%
%
%

\begin{tikzpicture}[scale=0.6]

\pgfdeclarelayer{bg}    
\pgfdeclarelayer{edges}
\pgfdeclarelayer{points}
\pgfsetlayers{bg,edges,points}

\def\rows{4}
\def\xscale{0.6}
\pgfmathtruncatemacro{\rowsminusone}{\rows-1}
\pgfmathtruncatemacro{\columns}{2^\rows/2}

\foreach \r in {0,1,...,\rows}{
	\begin{pgfonlayer}{bg}
	\draw (0,\r) -- (\xscale*\columns, \r);
	\end{pgfonlayer}
	\ifthenelse{\r < \rows}{
		\pgfmathtruncatemacro{\z}{\columns/2^\r}
		\foreach \i in {0,1,...,\z}{
			\begin{pgfonlayer}{bg}
			\draw (\xscale*2^\r * \i, \r) -- (\xscale*2^\r * \i, \r+1);
			\end{pgfonlayer}
		}
	}{}
}

\foreach \x/\y in {1/0, 2/0, 3/0, 4/0, 0/1, 1/1, 2/1, 0/2, 1/2}{
	\begin{pgfonlayer}{points}
			\draw[blue,fill=white, line width=0.9pt] (\xscale*2^\y * \x + 0.5*\xscale*2^\y, \y+ 0.5) circle [radius=0.08cm];
	\end{pgfonlayer}
	\begin{pgfonlayer}{edges}
			\draw[blue,line width=0.9pt] (\xscale*2^\y * \x + 0.5*\xscale*2^\y, \y+ 0.5) -- (\xscale*3,1.5);
	\end{pgfonlayer}
}

\end{tikzpicture}
\qquad
\begin{tikzpicture}

\pgfdeclarelayer{bg}    
\pgfdeclarelayer{edges}
\pgfdeclarelayer{points}
\pgfsetlayers{bg,edges,points}

\def\rows{4}
\def\xscale{0.6}
\pgfmathtruncatemacro{\rowsminusone}{\rows-1}
\pgfmathtruncatemacro{\columns}{2^\rows/2}

\foreach \r in {0,1,...,\rows}{
	\begin{pgfonlayer}{bg}
	\draw (0,\r) -- (\xscale*\columns, \r);
	\end{pgfonlayer}
	\ifthenelse{\r < \rows}{
		\pgfmathtruncatemacro{\z}{\columns/2^\r}
		\foreach \i in {0,1,...,\z}{
			\begin{pgfonlayer}{bg}
			\draw (\xscale*2^\r * \i, \r) -- (\xscale*2^\r * \i, \r+1);
			\end{pgfonlayer}
		}
	}{}
}

\foreach \x/\y in {2/0, 3/0, 0/1, 1/1, 2/1, 0/2}{
	\begin{pgfonlayer}{points}
			\draw[blue,fill=white, line width=0.9pt] (\xscale*2^\y * \x + 0.5*\xscale*2^\y, \y+ 0.5) circle [radius=0.08cm];
	\end{pgfonlayer}
	\begin{pgfonlayer}{edges}
			\draw[blue,line width=0.9pt] (\xscale*2^\y * \x + 0.5*\xscale*2^\y, \y+ 0.5) -- (\xscale*3,1.5);
	\end{pgfonlayer}
}

\end{tikzpicture}
\caption{The connection rules of $\mathcal{B}$ and $\mathcal{B}^\ast$. Left: a box with its $8$ neighbors in $\mathcal{B}$. 
Right: a box with its $5$ neighbors in $\mathcal{B}^\ast$.\label{fig:neighbors}}
\end{figure}

The boxes are the vertices of a graph $\mathcal{B}$\index{$\mathcal{B}$} in which two boxes are connected\index{box!connection rule} 
if they share at least a corner (Figure \ref{fig:neighbors}, left). Here we identify the left and right edge of $\mathcal{E}_R$ with 
each other, so that (for example) also the leftmost and rightmost box in each layer become neighbors. 
The dissection has the following properties:

\begin{lem} \label{lem:B-properties}
The following hold for $\Bcal$ and $\Gamma$:
\begin{enumerate}[(i)]
\item\label{lem:item:B-connections} If vertices of $\Gamma$ lie in boxes that are neighbors in $\mathcal{B}$, then they are 
connected by an edge in $\Gamma$. 
\item\label{lem:item:B-upperhalf} The number of boxes that lie (partly) above the line $y=R/2$ is at most $63$.
\end{enumerate}
\end{lem}

\begin{proof}
We start with \eqref{lem:item:B-connections}. 
Consider two points $(x,y)$ and $(x',y')$ that lie in boxes that are neighbors in $\mathcal{B}$. Suppose that the lowest of these two points lies 
in $L_i$. Then $y,y'\geq i\log(2)$. Furthermore, the horizontal distance between $(x,y)$ and $(x',y')$ is at most $3$ times the width of a 
box in $L_i$. It follows that
\[
|x-x'|_{\pi e^{R/2}} \leq 3 \cdot 2^i \cdot b \leq 2^i \leq e^{\frac12(y+y')},
\]
so $(x,y)$ and $(x',y')$ are indeed connected in $\Gamma$. 

To show \eqref{lem:item:B-upperhalf}, we note that the first layer $L_i$ that extends above the line $y=R/2$ has index 
$i = \lfloor \frac{R/2}{\log 2} \rfloor$. Therefore, we must count the boxes in the layers $L_i$, $L_{i+1}$, \ldots, $L_{\ell}$, of which 
there are $2^{\ell-i+1} - 1$. We have
\[
\ell-i+1 = \left\lfloor \frac{\log(6\pi) + R/2}{\log 2} \right\rfloor 
- \left\lfloor \frac{R/2}{\log 2} \right\rfloor + 1 \leq \left\lceil \frac{\log(6\pi)}{\log 2} \right\rceil + 1 = 6,
\]
so there are indeed at most $2^6-1=63$ boxes that extend above the line $y=R/2$.
\end{proof}

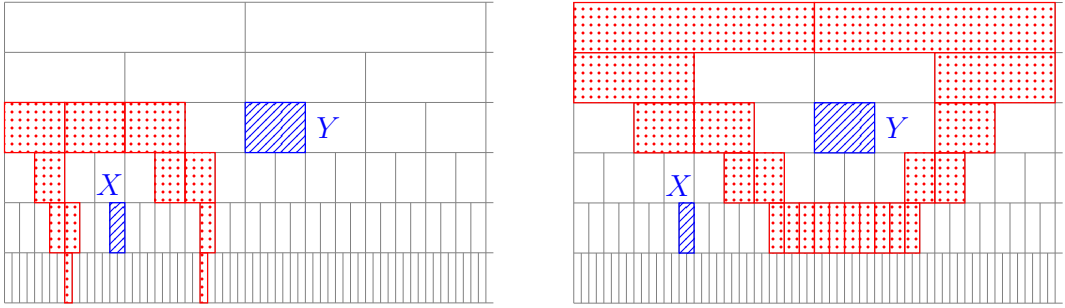
\begin{figure}[h!] \centering
\resizebox{0.4\linewidth}{!}{
\begin{tikzpicture}

\pgfmathsetseed{1138}

\pgfdeclarelayer{bg}    
\pgfdeclarelayer{edges}
\pgfdeclarelayer{text}
\pgfsetlayers{bg,main,edges,text}

\def\xscale{0.15}
\def\columns{65}
\def\rows{6}
\def\density{0.5}
\def\rowsh{4}
\pgfmathtruncatemacro{\rowsminusone}{\rows - 1}
\pgfmathtruncatemacro{\columnsa}{\columns/2}

\tikzset{%
	examplebox/.style = {color=blue,line width=0.7pt,pattern=north east lines,pattern color=blue},
    neighbor/.style={color=red,line width=0.7pt,pattern = dots,pattern color=red},
    gridstyle/.style={color=gray}
}

\foreach \i/\j in {1/7, 3/4}{
\begin{pgfonlayer}{edges}
\pgfmathsetmacro{\maxx}{min(\xscale*\j*2^\i+\xscale*2^\i,\columns*\xscale)}
\draw[examplebox] (\xscale*\j*2^\i,\i) rectangle (\maxx,\i + 1);
\end{pgfonlayer}
}

\node [above, blue] at (14*\xscale, 2) {\Large $X$};
\node [right, blue] at (40.5*\xscale, 3.5) {\Large $Y$};

\foreach \i/\j in {0/8, 1/4, 1/3, 2/1, 3/0, 3/1, 3/2, 2/5, 2/6, 1/13, 0/26}{
\begin{pgfonlayer}{edges}
\pgfmathsetmacro{\maxx}{min(\xscale*\j*2^\i+\xscale*2^\i,\columns*\xscale)}
\draw[neighbor] (\xscale*\j*2^\i,\i) rectangle (\maxx,\i + 1);
\end{pgfonlayer}
}

%

\foreach \x in {0,...,\rows}{
\draw [gridstyle] (0, \x) -- (\columns*\xscale, \x);
}


\foreach \y in {0,...,\rowsminusone} {
    \pgfmathtruncatemacro{\columnsa}{\columns/2^\y}
	\foreach \x in {0,...,\columnsa}{
	\pgfmathsetmacro{\z}{2^\y*\x} 
	\pgfmathtruncatemacro{\w}{\z}
	\ifthenelse{\w < \columns}{
	\draw [gridstyle] (\z*\xscale,\y) -- (\z*\xscale,\y+1);
	}{}
	} 
}

\end{tikzpicture}}$\qquad$
\resizebox{0.4\linewidth}{!}{
\begin{tikzpicture}

\pgfmathsetseed{1138}

\pgfdeclarelayer{bg}    
\pgfdeclarelayer{edges}
\pgfdeclarelayer{text}
\pgfsetlayers{bg,main,edges,text}

\def\xscale{0.15}
\def\columns{65}
\def\rows{6}
\def\density{0.5}
\def\rowsh{4}
\pgfmathtruncatemacro{\rowsminusone}{\rows - 1}
\pgfmathtruncatemacro{\columnsa}{\columns/2}

\node [above, blue] at (14*\xscale, 2) {\Large $X$};
\node [right, blue] at (40.5*\xscale, 3.5) {\Large $Y$};

\tikzset{%
	examplebox/.style = {color=blue,line width=0.7pt,pattern=north east lines,pattern color=blue},
    neighbor/.style={color=red,line width=0.7pt,pattern = dots,pattern color=red},
    gridstyle/.style={color=gray}
}

\foreach \i/\j in {1/7, 3/4}{
\begin{pgfonlayer}{edges}
\pgfmathsetmacro{\maxx}{min(\xscale*\j*2^\i+\xscale*2^\i,\columns*\xscale)}
\draw[examplebox] (\xscale*\j*2^\i,\i) rectangle (\maxx,\i + 1);
\end{pgfonlayer}
}

\foreach \i/\j in {3/2, 3/1, 4/0, 5/0, 5/1, 4/3, 3/6, 2/12, 2/11, 1/22, 1/21, 1/20, 1/19, 1/18, 1/17, 1/16, 1/15, 1/14, 1/13, 2/6, 2/5}{
\begin{pgfonlayer}{edges}
\pgfmathsetmacro{\maxx}{min(\xscale*\j*2^\i+\xscale*2^\i,\columns*\xscale)}
\draw[neighbor] (\xscale*\j*2^\i,\i) rectangle (\maxx,\i + 1);
\end{pgfonlayer}
}

%

\foreach \x in {0,...,\rows}{
\draw [gridstyle] (0, \x) -- (\columns*\xscale, \x);
}


\foreach \y in {0,...,\rowsminusone} {
    \pgfmathtruncatemacro{\columnsa}{\columns/2^\y}
	\foreach \x in {0,...,\columnsa}{
	\pgfmathsetmacro{\z}{2^\y*\x} 
	\pgfmathtruncatemacro{\w}{\z}
	\ifthenelse{\w < \columns}{
	\draw [gridstyle] (\z*\xscale,\y) -- (\z*\xscale,\y+1);
	}{}
	} 
}

\end{tikzpicture}}
\caption{If two blue boxes $X$ and $Y$ are not connected by a path of blue (striped) boxes, then a red (dotted) walk exists that 
intersects every path in $\mathcal{B}$ from $X$ from $Y$ (Lemma \ref{obs:barrier}). 
This walk can be chosen such that it either connects two boxes in $L_0$ (left) or is cyclic (right). \label{fig:obstruction}}
\end{figure}

Let $\mathcal{B}^\ast$ be the subgraph of $\mathcal{B}$ where we remove the edges between boxes that have only a single point in 
common (Figure \ref{fig:neighbors}, right). 
Note that $\mathcal{B}^\ast$ is a planar graph and that $\mathcal{B}$ is obtained from $\mathcal{B}^\ast$ by adding the diagonals of 
each face (\cite{Kesten} deals with a more general notation of \emph{matching pairs} of graphs). 
We make the following observation (Figure \ref{fig:obstruction}; compare Proposition 2.1 in \cite{Kesten}) for later reference.

\begin{lem}\label{obs:barrier}
Suppose each box in $\mathcal{B}$ is colored red or blue. 
If there is no path of blue boxes in $\mathcal{B}$ between two blue boxes $X$ and $Y$, then $\mathcal{B}^\ast$ 
(and hence also $\mathcal{B}$) contains a walk of red boxes $Q$ that intersects every walk (and hence also every path) 
in $\mathcal{B}$ from $X$ to $Y$. \hfill \blokje
\end{lem}

\noindent
We leave the straightforward proof of this last lemma to the reader. It can for instance be derived quite succinctly from the
Jordan curve theorem. A proof can be found in the MSc thesis of the second author~\cite{MerlijnThesis}.

\subsection{Constructing short paths\label{sec:paths}}

We will use the boxes defined in the previous subsection to construct short paths between vertices of $\Gamma$. 
Recall that $V \subset \mathcal{E}_R$ is an arbitrary finite set of points and $\Gamma$ is the graph with vertex set $V$ and 
connection rule $(x,y) \sim (x',y') \iff |x-x'|_{\pi e^{R/2}} \leq e^{\frac12(y+y')}$. We will also make use of the dissection into boxes introduced in the previous section. 
A box is called \emph{active} if it contains at least one vertex of $\Gamma$ and \emph{inactive} otherwise.

\begin{figure}[h!]\centering
\resizebox{!}{0.4\linewidth}{
\begin{tikzpicture}

\pgfdeclarelayer{bg}    
\pgfdeclarelayer{edges}
\pgfdeclarelayer{text}
\pgfsetlayers{bg,edges,main,text}

\def\xscale{0.07}
\def\columns{135}
\def\rows{9}
\def\density{0.5}
\def\rowsh{4}
\pgfmathtruncatemacro{\rowsminusone}{\rows - 1}
\pgfmathtruncatemacro{\rowsplusone}{\rows + 1}
\pgfmathtruncatemacro{\rowsplustwo}{\rows + 2}
\pgfmathtruncatemacro{\columnsa}{\columns/2}


\begin{pgfonlayer}{bg}
\foreach \x in {0,...,\rows}{
\draw [white] (0, \x) -- (\columns*\xscale, \x);
}

\foreach \y in {0,...,\rowsminusone} {
    \pgfmathtruncatemacro{\columnsa}{\columns/2^\y}
	\foreach \x in {0,...,\columnsa}{
	\pgfmathsetmacro{\z}{2^\y*\x} 
	\pgfmathtruncatemacro{\w}{\z}
	\ifthenelse{\w < \columns}{
	\draw [white] (\z*\xscale,\y) -- (\z*\xscale,\y+1);
	}{}
	} 
}
\end{pgfonlayer}
%

\tikzset{%
	examplebox/.style = {color=gray!30!white},
    neighbor/.style={color=black, pattern=north east lines,pattern color=gray},
}

\begin{pgfonlayer}{text}
\node [black, right] at (50*\xscale,8.5) {$L(A,A')$};
\draw [black, line width=0.8pt] (82*\xscale,7.5) -- (75*\xscale,8.5);
\end{pgfonlayer}

\foreach \i/\j in {2/10, 3/5, 4/2, 5/1, 6/0, 1/55, 2/27, 3/13, 4/6, 5/3, 6/1, 7/0}{
\begin{pgfonlayer}{bg}
\pgfmathsetmacro{\maxx}{min(\xscale*\j*2^\i+\xscale*2^\i,\columns*\xscale)}
\draw[neighbor] (\xscale*\j*2^\i,\i) rectangle (\maxx,\i + 1);
\end{pgfonlayer}
}

\begin{pgfonlayer}{text}
\node [below, black] at (37*\xscale,2.7) {$A$};
\node [right, black] at (112*\xscale,1.5) {$A'$};
\end{pgfonlayer}


%
%
%
%
%

\end{tikzpicture}}$\qquad$\resizebox{!}{0.4\linewidth}{
\begin{tikzpicture}

\pgfdeclarelayer{bg}    
\pgfdeclarelayer{edges}
\pgfdeclarelayer{text}
\pgfsetlayers{bg,edges,main,text}

\def\xscale{0.07}
\def\columns{135}
\def\rows{9}
\def\density{0.5}
\def\rowsh{4}
\pgfmathtruncatemacro{\rowsminusone}{\rows - 1}
\pgfmathtruncatemacro{\rowsplusone}{\rows + 1}
\pgfmathtruncatemacro{\rowsplustwo}{\rows + 2}
\pgfmathtruncatemacro{\columnsa}{\columns/2}

\tikzset{%
	examplebox/.style = {color=gray!30!white},
    neighbor/.style={pattern=north east lines,pattern color=gray},
    inactive/.style={color=white,fill=white},
    w/.style={pattern=none}
}

\begin{pgfonlayer}{bg}
\fill [black!20!white] (-0.25,0) rectangle (\columns*\xscale,\rows);
\end{pgfonlayer}


%

\foreach \i/\j in {3/5,4/2,4/3,3/4,3/7,3/8,2/16,2/17,3/3,2/7,2/6,1/14,1/13,2/5,2/4,1/7,1/8,5/0,1/34,1/33,1/32,4/6,5/3,3/11,2/22,2/27,2/28,2/29,2/30,3/15,7/0,
1/61,1/62,1/63,1/64}{
\begin{pgfonlayer}{main}
\pgfmathsetmacro{\maxx}{min(\xscale*\j*2^\i+\xscale*2^\i,\columns*\xscale)}
\draw[inactive] (\xscale*\j*2^\i,\i) rectangle (\maxx,\i + 1);
\end{pgfonlayer}
}

\foreach \i/\j in {2/20,1/40,2/33,1/66,0/3,0/4,0/5,0/6,0/7,0/40,0/41,0/42,0/44,0/45,0/46,0/47,0/48,
0/50,0/51,0/54,0/55,0/56,0/80,0/81,0/82,0/83,0/86,0/87,0/88,0/89,0/90,0/95,0/96,0/97,0/100,0/101,0/102,
0/116,0/117,0/118,0/115,
3/0,2/0,2/1,1/3,1/4,1/24,1/25,1/26,1/27,1/47,1/48,1/49,1/50,1/37,1/38,2/19, 0/15,0/16,0/17,0/18,0/19,0/62,0/63,0/64,0/65,0/66}{
\begin{pgfonlayer}{main}
\pgfmathsetmacro{\maxx}{min(\xscale*\j*2^\i+\xscale*2^\i,\columns*\xscale)}
\draw[inactive] (\xscale*\j*2^\i,\i) rectangle (\maxx,\i + 1);
\end{pgfonlayer}
}

\foreach \i/\j in {2/10, 3/5, 4/2, 5/1, 6/0, 1/55, 2/27, 3/13, 4/6, 5/3, 6/1, 7/0}{
\begin{pgfonlayer}{main}
\pgfmathsetmacro{\maxx}{min(\xscale*\j*2^\i+\xscale*2^\i,\columns*\xscale)}
\xdef\d{0}
\fill [neighbor] (\xscale*\j*2^\i-\d,\i-\d) rectangle (\maxx+\d,\i + 1+\d);
\end{pgfonlayer}
}

\begin{pgfonlayer}{text}
\node [below, black] at (37*\xscale,2.7) {$A$};
\node [right, black] at (111*\xscale,1.5) {$A'$};
\end{pgfonlayer}


\xdef\xs{\xscale}
\draw[black, line width=1.2pt] (15*\xs,0) -- (15*\xs,1) -- (14*\xs,1)--(14*\xs,2) -- (16*\xs,2) -- (16*\xs,3) -- (24*\xs,3)--(24*\xs,4)--(32*\xs,4)--(32*\xs,5)--(0,5)--(0,8)--(128*\xs,8) -- (128*\xs,5) -- (112*\xs,5)--(112*\xs,3) -- (120*\xs,3)--(120*\xs,4) -- (128*\xs,4) -- (128*\xs,3) -- (124*\xs,3) -- (124*\xs,2) -- (130*\xs,2) -- (130*\xs,1) -- (122*\xs,1) -- (122*\xs,2) -- (112*\xs,2)--(112*\xs,1) -- (110*\xs,1)--(110*\xs,2)--(108*\xs,2) -- (108*\xs,3) -- (104*\xs,3) -- (104*\xs,4) -- (96*\xs,4) -- (96*\xs,3) -- (92*\xs,3) -- (92*\xs,2) -- (88*\xs,2) -- (88*\xs,4) -- (96*\xs,4) -- (96*\xs,6) -- (64*\xs,6)--(64*\xs,4) -- (72*\xs,4)--(72*\xs,2)--(70*\xs,2) -- (70*\xs,1)--(67*\xs,1)--(67*\xs,0)--(62*\xs,0)--(62*\xs,1)--(64*\xs,1)--(64*\xs,3)--(56*\xs,3)--(56*\xs,4) -- (48*\xs,4)--(48*\xs,3) -- (44*\xs,3)--(44*\xs,2)--(40*\xs,2)--(40*\xs,3) -- (32*\xs,3)--(32*\xs,2)--(30*\xs,2)--(30*\xs,1)--(26*\xs,1)--(26*\xs,2)--(18*\xs,2)--(18*\xs,1)--(20*\xs,1) --(20*\xs,0)--(15*\xs,0);

\begin{pgfonlayer}{text}
\node [black, right] at (48*\xscale,8.5) {$W(A,A')$};
\draw [black, line width=0.8pt] (78.5*\xscale,8) -- (75*\xscale,8.5);
\end{pgfonlayer}

\end{tikzpicture}}
\caption{Two boxes $A, A'$ in $\mathcal{B}$ and the path $L(A,A')$ that connects them (left). 
We can form $L(A,A')$ by concatenating the shortest paths from $A$ and $A'$ to the lowest box lying above both $A$ and $A'$. 
In the right image active boxes are colored gray and inactive boxes are colored white. The union of $L(A,A')$ and the inactive components 
intersecting $L(A,A')$ is called $W(A,A')$ and outlined in black. \label{fig:LW}}
\end{figure}
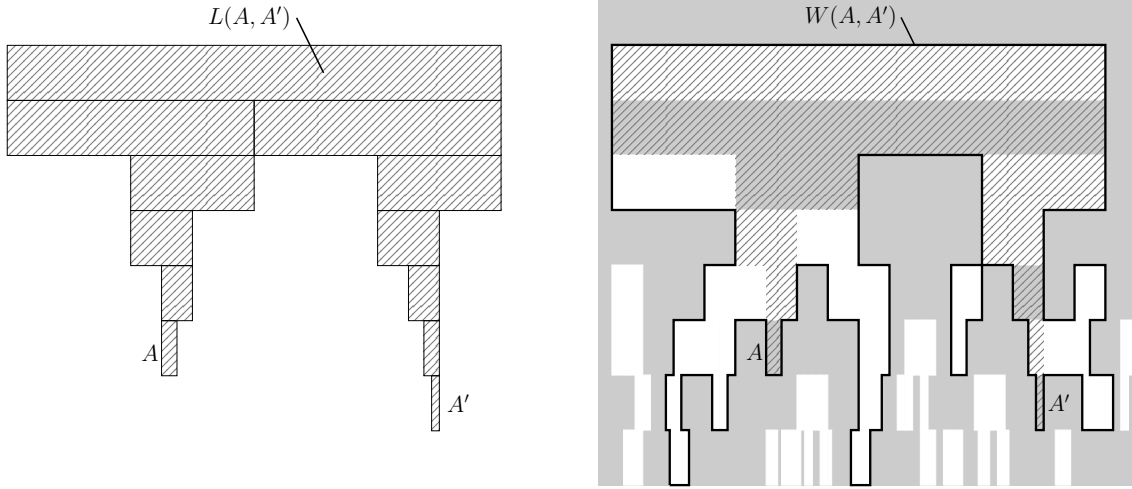

Suppose $x$ and $x'$ are two vertices of $\Gamma$ that lie in the same component. How can we find a short path from $x$ to $x'$? A natural strategy would be to follow a short path of boxes from the box $A$ containing $x$ to the box $A'$ containing $x'$. These boxes are connected by a path $L(A,A')$ of length at most $2R$ (Figure \ref{fig:LW}, left). 
If all the boxes in $L(A,A')$ are active,  Lemma \ref{lem:B-properties}\eqref{lem:item:B-connections} immediately yields a path in $\Gamma$ from $x$ to $x'$ of length at most $2R$, which is a path of the desired length. The situation is more difficult if we also encounter inactive boxes, and modifying the path to avoid inactive boxes may be impossible because a path of active boxes connecting $A$ to $A'$ may fail to exist. Nevertheless, it turns out that the graph-theoretic distance between $x$ and $x'$ can be bounded in terms of the size of inactive regions one encounters when following $L(A,A')$. 

To make this precise, we define $W(A,A')$ to be the set of boxes that either lie in $L(A,A')$ or from which an inactive path (i.e.\ a path of inactive boxes) exists to 
a box in $L(A,A')$ (Figure \ref{fig:LW}, right). Note that $W(A,A')$ is a connected subset of $\mathcal{B}$, consisting of all boxes 
in $L(A,A')$ and all inactive components intersecting $W(A,A')$ (by an inactive component\index{inactive component} we mean a component 
of the induced subgraph of $\mathcal{B}$ on the inactive boxes). The main result of this section is that the graph-theoretic distance between $x$ and $x'$ is bounded by the size of $W(A,A')$.

Before we continue, we first recall some geometric properties of the graph $\Gamma$:\index{idealized model!geometric properties}

\begin{lem}[\cite{FM}, Lemma 3]\label{lem:geometry}
Let $x$, $y$, $z$, $w \in V$. 
\begin{itemize}
\item[(i)] If $xy \in E(\Gamma)$ and $z$ lies above the line segment $[x,y]$ 
(i.e. $[x,y]$ intersects the segment joining $z$ and the projection of $z$ onto the horizontal axis), then at least one 
of $xz$ and $yz$ is also present in $\Gamma$.
\item[(ii)] If $xy$, $wz \in E(\Gamma)$ and the segments $[x,y]$ and $[z,w]$ intersect, then at least one of the
edges $xw$, $xz$, $yw$ and $yz$ is also present in $\Gamma$. In particular, $\{x,y,z,w\}$ is a connected subset of $\Gamma$.
\end{itemize}
\end{lem}

We now prove a lemma that allows us to compare paths in $\Gamma$ with walks in $\mathcal{B}$. 
This will enable us to translate information about $\Gamma$ (such as that two boxes contain vertices in the same component of $\Gamma$) to information about the states of the boxes. 

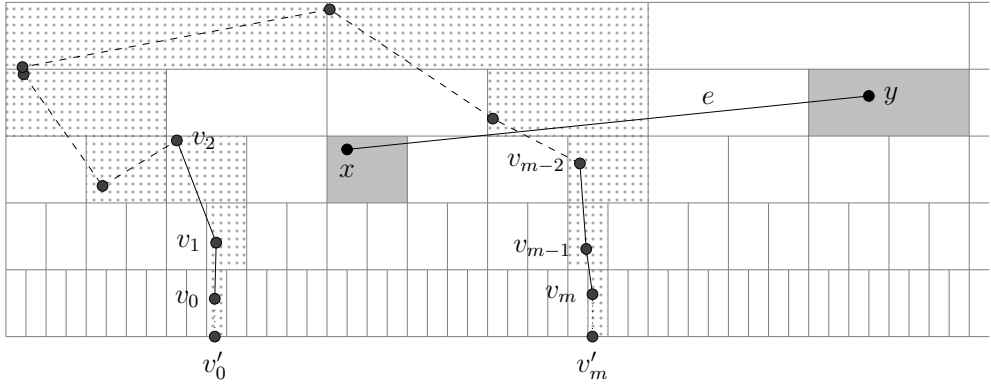
\begin{figure}[h!] \centering
\resizebox{0.8\linewidth}{!}{

\begin{tikzpicture}

\pgfdeclarelayer{bg}    
\pgfdeclarelayer{edges}
\pgfdeclarelayer{text}
\pgfsetlayers{bg,edges,main,text}

\def\xscale{0.3}
\def\columns{49}
\def\rows{5}
\def\n{8}
\pgfmathtruncatemacro{\rowsminusone}{\rows - 1}
\pgfmathtruncatemacro{\rowsplusone}{\rows + 1}
\pgfmathtruncatemacro{\rowsplustwo}{\rows + 2}
\pgfmathtruncatemacro{\columnsa}{\columns/2}


{

\tikzset{%
	examplebox/.style = {color=gray!50!white},
    neighbor/.style={pattern color=gray!70!white,pattern=dots},
}

\begin{pgfonlayer}{bg}
\fill[examplebox] (16*\xscale,2) rectangle (20*\xscale,3);
\fill[examplebox] (40*\xscale,3) rectangle (48*\xscale,4);
\end{pgfonlayer}

\draw[black,fill=black] (17*\xscale,2.8) circle [radius=0.08]; 
\draw[below=0.08cm] (17*\xscale,2.8) node {$x$};
\draw[black,fill=black] (43*\xscale,3.6) circle [radius=0.08]; 
\draw[right=0.08cm] (43*\xscale,3.6) node {$y$};
\draw[black] (17*\xscale,2.8) -- (43*\xscale,3.6);
\draw[above=0.05cm] (35*\xscale,3.29) node {$e$};

\xdef\k{0}

\foreach \i/\j/\k in {0/10/0, 1/5/1, 2/2/2, 2/1/3, 3/0/4, 4/0/5, 4/1/6, 3/3/7, 2/7/8, 1/14/9, 0/29/10}{
\begin{pgfonlayer}{bg}
\fill[neighbor] (\xscale*\j*2^\i,\i) rectangle (\xscale*\j*2^\i+\xscale*2^\i,\i + 1);
\end{pgfonlayer}
\pgfmathsetmacro{\randomx}{rnd}
\pgfmathsetmacro{\randomy}{rnd}
\ifthenelse{\k=0 \OR \k=6}{\pgfmathsetmacro{\randomy}{0.5+0.5*rnd}}{}
\ifthenelse{\k=7 \OR \k=8}{\pgfmathsetmacro{\randomy}{0.9*rnd}}{}
\ifthenelse{\k=10}{\pgfmathsetmacro{\randomy}{0.5+0.5*rnd}}{}
\xdef\newx{\xscale*\j*2^\i+\randomx*\xscale}
\xdef\newy{\i+\randomy}

\begin{pgfonlayer}{edges}
\ifthenelse{\k>0}{
\ifthenelse{\k>2 \AND \k<9}{
\draw[black,dashed] (\newx,\newy) -- (\oldx,\oldy);}{
\draw[black] (\newx,\newy) -- (\oldx,\oldy);}
}{
\draw[black,dotted] (\newx, \newy) -- (\newx, 0);
\draw [black,fill=black!75!white] (\newx,0) circle [radius=0.08];
\begin{pgfonlayer}{text}
\draw [below=0.08cm] (\newx,0) node {$v_0'$};
\end{pgfonlayer}
}
\end{pgfonlayer}
\begin{pgfonlayer}{main}
\xdef\oldx{\newx}
\xdef\oldy{\newy}
\begin{pgfonlayer}{text}
\draw [black,fill=black!75!white] (\newx,\newy) circle [radius=0.08];
\end{pgfonlayer}
\ifthenelse{\k<3}{
\begin{pgfonlayer}{text}
\ifthenelse{\k=2}{
\draw [right=0.08cm] (\newx,\newy) node {$v_{\k}$};
}{
\draw [left=0.08cm] (\newx,\newy) node {$v_{\k}$};}
\end{pgfonlayer}}{}
\ifthenelse{\k>7}{
	\pgfmathtruncatemacro{\l}{10-\k}
	\ifthenelse{\k=7}{
	\begin{pgfonlayer}{text}
		\draw [left=0.08cm] (\newx,\newy) node {$v_{m-\l}$};
		\end{pgfonlayer}
	}{
	\ifthenelse{\k>9}{
	\begin{pgfonlayer}{text} \draw [left=0.08cm] (\newx,\newy) node {$v_{m}$}; \end{pgfonlayer}}{
	\draw [left=0.08cm] (\newx,\newy) node {$v_{m-\l}$};}}}{}
\end{pgfonlayer}
}

\draw[black,dotted] (\newx, \newy) -- (\newx, 0);
\begin{pgfonlayer}{text}
\draw [black,fill=black!75!white] (\newx,0) circle [radius=0.08];
\draw [below=0.08cm] (\newx,0) node {$v_m'$};
\end{pgfonlayer}
}

%
%
%
%


\begin{pgfonlayer}{bg}

\foreach \x in {0,...,\rows}{
\draw [gray] (0, \x) -- (\columns*\xscale, \x);
}

%


\foreach \y in {0,...,\rowsminusone} {
    \pgfmathtruncatemacro{\columnsa}{\columns/2^\y}
	\foreach \x in {0,...,\columnsa}{
	\pgfmathsetmacro{\z}{2^\y*\x}
	\pgfmathtruncatemacro{\w}{\z}
	\ifthenelse{\w < \columns}{
	\draw [gray] (\z*\xscale,\y) -- (\z*\xscale,\y+1);
	}{}
	} 
}

\end{pgfonlayer}

%

%
%

\def\bottomy{0}

%
%

\end{tikzpicture}}
\caption{Proof of Lemma \ref{lem:edge-replace}. The edge $e$ of $\Gamma$ connects vertices $x$ and $y$. If a red (dotted) walk of boxes 
exists that separates the boxes containing $x$ and $y$, then $e$ intersects one of the 
segments $[v_i,v_{i+1}]$, $[v_0,v_0']$ or $[v_m,v_m']$. 
This contradicts the assumption that no red box contains a neighbor of $x$ or $y$. \label{fig:lem:edge-replace}}
\end{figure}

\begin{lem}\label{lem:Qpath}\label{lem:edge-replace}
Suppose boxes $X,Y \in \mathcal{B}$ contain vertices $x,y \in V$ respectively that lie in the same component of $\Gamma$. 
Then $\mathcal{B}$ contains a walk $X= B_0$, $B_1$, \ldots, $B_n = Y$ with the following property: \begin{itemize}
\item[($\ast$)] if $B_i$ and $B_j$ are active but $B_{i+1}$, $B_{i+2}$, \ldots, $B_{j-1}$ are not, then $\Gamma$ has 
vertices $a \in B_i$, $b \in B_j$ that are connected in $\Gamma$ by a path of length at most $3$.
\end{itemize}
\end{lem}

\begin{proof}
We prove the statement by induction on the length of the shortest path from $x$ to $y$ in $\Gamma$. 

First suppose that this length is $1$, so that there is an edge $e$ connecting $x$ and $y$. We claim that a 
walk $X=B_0$, $B_1$, \ldots, $B_n = Y$ in $\mathcal{B}$ exists with the property that if $B_i$ is active, then $B_i$ contains a 
neighbor of $x$ or $y$. 
For this we use Lemma \ref{obs:barrier}. We color a box blue if it is either {\bf a)} inactive or {\bf b)} active and it contains 
a neighbor of $x$ or $y$. 
All other boxes are colored red. 
Note that $X$ and $Y$ are blue, because $X$ contains 
a neighbor of $y$ (namely $x$) and $Y$ contains a neighbor of $x$ (namely $y$). 
We intend to show that $\mathcal{B}$ contains a blue path 
from $X$ to $Y$. Aiming for a contradiction, we suppose that this is not the case. 
By Lemma \ref{obs:barrier}, there must then exist a red walk $S=S_0$, $S_1$, \ldots, $S_m$ that intersects each path 
in $\mathcal{B}$ from $X$ to $Y$. 
If we remove $S$ from $\mathcal{E}_R$ then $\mathcal{E}_R \backslash S$ falls apart in a number of components. Because there is no path in 
$\mathcal{B}$ from $X$ to $Y$ that does not intersect $S$, $X$ and $Y$ lie in different components. 
(We say $S$ \emph{separates} $X$ and $Y$.)
We choose vertices $v_i \in V \cap S_i$ for all $i$ 
(these vertices exist because all red  boxes are active; see Figure \ref{fig:lem:edge-replace}). 
By Lemma \ref{lem:B-properties}\eqref{lem:item:B-connections}, $v_i$ and $v_{i+1}$ are neighbors in $\Gamma$ for each $i$.

We may assume that either $S_0$ and $S_m$ are both boxes in the lowest layer $L_0$, or $S_0$ and $S_m$ are adjacent 
in $\mathcal{B}$ (Figure \ref{fig:obstruction}). 
In the latter case, we consider the polygonal curve $\gamma$ consisting of the line segments $[v_0,v_1]$, $[v_1,v_2]$, \ldots, $[v_m,v_0]$.
This polygonal curve consists of edges of $\Gamma$. 
Let us observe that each of these edges passes through boxes in $S$ and maybe also boxes adjacent to boxes in $S$, but the edges cannot intersect 
any box that is neither on $S$ nor adjacent to a box of $S$. 
So in particular, none of these edges can pass through the box $X$, because $X$ is not adjacent to a box in $S$ 
(this box should then have been blue by Lemma \ref{lem:B-properties}\eqref{lem:item:B-connections}). 
From this it follows that $\gamma$ also separates $x$ and $y$. 
Therefore, the edge $e$ crosses an edge $[v_i,v_{i+1}]$ of $\Gamma$ (Figure \ref{fig:lem:edge-replace}). 
By Lemma \ref{lem:geometry}(ii) this means that $v_i$ or $v_{i+1}$ neighbors $x$ or $y$, which is a contradiction because 
$v_i$ and $v_{i+1}$ do not lie in a blue box.

We are left with the case that $S_0$ and $S_m$ lie in the lowest layer $L_0$. 
Let $v_0'$ and $v_m'$ denote the projections of $v_0$ and $v_m$, respectively, on the horizontal axis. 
By an analogous argument, we find that the polygonal line through $v_0'$, $v_0$, $v_1$, \ldots, $v_m$, $v_m'$ separates $x$ and $y$. 
We now see that $e$ either crosses an edge $[v_i,v_{i+1}]$ (we then find a contradiction with Lemma \ref{lem:geometry}(ii)) or 
one of the segments $[v_0,v_0']$ and $[v_m,v_m']$ (we then find a contradiction with Lemma \ref{lem:geometry}(i)). 
From the contradiction we conclude that a blue path must exist connecting $X$ and $Y$.

We have now shown that if $x$ and $y$ are neighbors in $\Gamma$, there exists a walk $X = B_0$, $B_1$, \ldots, $B_n = Y$ such that 
the $B_i$ that are active contain a neighbor of $x$ or $y$. This means that if $B_i$, $B_{i+1}$, \ldots, $B_j$ are such that $B_i$ and $B_j$ are 
active but $B_{i+1}$, \ldots, $B_{j-1}$ are not, then $B_i$ contains a vertex $a$ that neighbors $x$ or $y$ and $B_j$ contains a vertex $b$ 
that neighbors $x$ or $y$. Now $d_{\Gamma}(a,b) \leq 3$ follows from the fact that both $a$ and $b$ neighbor an endpoint of the same edge $e$.
We conclude that if $x$ and $y$ are neighbors in $\Gamma$, then a walk satisfying ($\ast$) exists.

Now suppose that the statement holds 
whenever $x$ and $y$ satisfy $d_{\Gamma}(x,y) \leq k$ and consider two vertices $x$ and $y$ with $d_{\Gamma}(x,y) = k+1$. 
Choose a neighbor $y'$ of $y$ such that $d_{\Gamma}(x,y')=k$. Let $Y'$ be the active box containing $y'$. 
By the induction hypothesis, there exists walks from $X$ to $Y'$ and from $Y'$ to $Y$ satisfying ($\ast$). 
By concatenating these two walks we obtain a walk from $X$ to $Y$ satisfying ($\ast$), as desired.
\end{proof}

Note that by itself this lemma is insufficient to construct short paths, as the proof is non-constructive and there is no control over the number of boxes in the walk obtained. Nevertheless, we can use Lemma \ref{lem:edge-replace} to prove the main result of this section.

\begin{lem}\label{lem:distance-bound-1}
There exists a constant $c$ such that the following holds (for all finite $V \subseteq \Ecal_R$ with $\Gamma$ constructed as above).
If the vertices $x, x' \in V$ and the boxes $A, A' \in \mathcal{B}$ are such that
\begin{enumerate}
\item $x \in A$, $x' \in A'$, and; 
\item $x,x'$ lie in the same component of $\Gamma$, 
\end{enumerate}
then $d_\Gamma(x,x') \leq c |W(A,A')|$.
\end{lem}

\begin{proof}
We claim that there is a walk $S=S_0$, $S_1$, \ldots, $S_n$ in $\mathcal{B}$ from $A$ to $A'$ satisfying

\begin{itemize}
\item[(i)] if $S_i$ and $S_j$ are active but $S_{i+1}$, \ldots, $S_{j-1}$ are not, then $\Gamma$ has vertices $a \in S_i$, $b \in S_j$ that 
are connected in $\Gamma$ by a path of length at most $3$;
\item[(ii)] if $S_i$ is active, then either $S_i$ itself or an inactive box adjacent to $S_i$ belongs to $W(A,A')$.
\end{itemize}

\noindent
We define $\mathcal{B}_x$ to be the set of active boxes that contain vertices of the component of $\Gamma$ that contains $x$ and $x'$. 
By assumption we have $A, A' \in \mathcal{B}_x$. If $A$ and $A'$ are adjacent the existence of a walk $S$ satisfying (i) and (ii) is trivial, so we assume $A$ and $A'$ are not adjacent. The proof consists of proving the result for the case that $A$ and $A'$ are the only boxes in $L(A,A')$ that belong to $\mathcal{B}_x$, and then a straightforward extension to the general case.

If $A$ and $A'$ are the only boxes in $L(A,A')$ that belong to $\mathcal{B}_x$, then the boxes in between $A$ and $A'$ 
on $L(A, A')$ are either inactive, or they are active but contain vertices of a different component of $\Gamma$. 
Therefore, the box $B$ in $L(A,A')$ directly following $A$ must be inactive and belongs to some inactive component $F$ (recall that an inactive component 
is a component of the induced 
subgraph of $\mathcal{B}$ on the inactive boxes). We will prove the stronger statement that a walk $S=S_0$, $S_1$, \ldots, $S_n$ from $A$ to $A'$ 
exists satisfying (i) and

\begin{itemize}
\item[(ii')] if $S_i$ is active, then $S_i$ is adjacent to a box in $F$.
\end{itemize}

\noindent
By Lemma \ref{lem:Qpath} there exists a walk $S$ from $A$ to $A'$ satisfying (i). We will modify $S$ such that also (ii') holds. 
We proceed in two steps. In Step 1 we remove all inactive boxes in $S$ that are not in $F$. 
In Step 2 we remove all active boxes from $S$ that are not adjacent to a box in $F$. \vspace{0.6\baselineskip}

\textit{Step 1.} There is a walk $S$ satisfying (i) that contains no inactive boxes outside $F$. \vspace{0.6\baselineskip}

We start with the walk $S$ that Lemma \ref{lem:Qpath} provides. This walk satisfies (i). 
Suppose $S$ contains some inactive box $X$ not in $F$ (Figure \ref{fig:lemma8}, left). 
Because $B \in F$, there can then be no inactive path in $\mathcal{B}$ from $X$ to $B$. 
It follows from Lemma \ref{obs:barrier} that there is an active walk $Q$ that intersects all walks in $\mathcal{B}$ from $X$ to $B$ 
(we apply Lemma \ref{obs:barrier} with the inactive boxes colored blue and all other boxes colored red). 
One such walk from $X$ to $B$ is obtained by following $S$ towards $A$ (which is a neighbor of $B$). 
Another possible walk is obtained by first following $S$ towards $A'$ and then following $L(A,A')$ towards $B$. 
We define boxes $E$ and $E'$ such that $Q$ intersects the walk in $\mathcal{B}$ from $X$ to $B$ via $S$ and $A$ in $E$ and the walk 
in $\mathcal{B}$ from $X$ to $B$ via $S$, $A'$ and $L(A,A')$ in $E'$ (Figure \ref{fig:lemma8}, left). 
Because $E$ belongs to $S$, $E$ also belongs to $\mathcal{B}_x$. It follows that $E'$ also belongs to $\mathcal{B}_x$, which implies 
that $E'$ lies in $S$ (the boxes in $L(A,A')$ between $A$ and $A'$ do not lie in $\mathcal{B}_x$ by assumption).  
We see that $Q$ contains two active boxes $E$ and $E'$ that lie on either side of $X$. Because $Q$ contains only active boxes, we can replace 
the part of $S$ from $E$ to $E'$ by a walk of active boxes from $E$ to $E'$. 
Doing so we find a walk that still satisfies (i) but from which the box $X$ is removed. 
By repeatedly applying this procedure, we remove all such boxes $X$ from $S$. The resulting walk satisfies (i) and contains no inactive boxes 
outside $F$. \vspace{0.6\baselineskip}

\begin{figure}[h]\centering
\resizebox{!}{0.35\linewidth}{
\begin{tikzpicture}

\pgfdeclarelayer{bg}    
\pgfdeclarelayer{edges}
\pgfdeclarelayer{text}
\pgfsetlayers{bg,edges,main,text}

\def\xscale{0.1}
\def\columns{135}
\def\rows{9}
\def\density{0.5}
\def\rowsh{4}
\pgfmathtruncatemacro{\rowsminusone}{\rows - 1}
\pgfmathtruncatemacro{\rowsplusone}{\rows + 1}
\pgfmathtruncatemacro{\rowsplustwo}{\rows + 2}
\pgfmathtruncatemacro{\columnsa}{\columns/2}


%

\tikzset{%
	examplebox/.style = {color=gray!30!white},
    neighbor/.style={color=black,fill=white},
    active/.style={color=black,fill=gray!50!white},
}

\begin{pgfonlayer}{text}
\node [black, right] at (70*\xscale,6.3) {$L(A,A')$};
\end{pgfonlayer}

\begin{pgfonlayer}{text}
\node [below, black] at (37*\xscale,2.7) {$A$};
\node [below, black] at (37*\xscale,3.7) {$B$};
\node [below, black] at (58*\xscale,2) {$E$};
\node [below, black] at (97*\xscale,1.5) {$E'$};
\node [black] at (68*\xscale, 1.1) {$X$};
\node [right, black] at (112*\xscale,1.5) {$A'$};
\node [right, black] at (112*\xscale,2.5) {$B'$};
\node [black] at (76*\xscale, -0.2) {$S$};
\node [black] at (85*\xscale,4.1) {$Q$};
\end{pgfonlayer}

\draw[black, line width=0.5pt] (42*\xscale, 4) -- (42*\xscale, 6) -- (111*\xscale, 6) -- (111*\xscale, 3);

\draw [black] plot [smooth, tension=2] coordinates { (42*\xscale,2) (42*\xscale,1) (65*\xscale,2) (80*\xscale,0) (90*\xscale,2) (110*\xscale,0) (111*\xscale, 1)};

\draw [gray, dotted, line width = 0.8pt] plot [smooth, tension =2] coordinates { (71*\xscale,2) (60*\xscale,3.5) (48*\xscale,3.5) };

\draw [black, dashed] plot [smooth, tension =2] coordinates { (58*\xscale,3) (80*\xscale,4) (101*\xscale,2) };

\foreach \i/\j in {2/10, 3/5, 1/55, 2/27, 1/35, 2/14, 1/50}{
\pgfmathsetmacro{\maxx}{min(\xscale*\j*2^\i+\xscale*2^\i,\columns*\xscale)}
\draw[neighbor] (\xscale*\j*2^\i,\i) rectangle (\maxx,\i + 1);
}

\foreach \i/\j in {2/10, 1/55, 2/14, 1/50}{
\pgfmathsetmacro{\maxx}{min(\xscale*\j*2^\i+\xscale*2^\i,\columns*\xscale)}
\draw[active] (\xscale*\j*2^\i,\i) rectangle (\maxx,\i + 1);
}

\end{tikzpicture}} \quad \resizebox{!}{0.35\linewidth}{
\begin{tikzpicture}

\pgfdeclarelayer{bg}    
\pgfdeclarelayer{edges}
\pgfdeclarelayer{text}
\pgfsetlayers{bg,edges,main,text}

\def\xscale{0.1}
\def\columns{135}
\def\rows{9}
\def\density{0.5}
\def\rowsh{4}
\pgfmathtruncatemacro{\rowsminusone}{\rows - 1}
\pgfmathtruncatemacro{\rowsplusone}{\rows + 1}
\pgfmathtruncatemacro{\rowsplustwo}{\rows + 2}
\pgfmathtruncatemacro{\columnsa}{\columns/2}


%

\tikzset{%
	examplebox/.style = {color=gray!30!white},
    neighbor/.style={color=black,fill=white},
    active/.style={color=black,fill=gray!50!white},
    F/.style={color=black,fill=gray!50!white,pattern=north east lines},
}

\begin{pgfonlayer}{text}
\node [black, right] at (70*\xscale,6.3) {$L(A,A')$};
\end{pgfonlayer}

\begin{pgfonlayer}{text}
\node [below, black] at (37*\xscale,2.7) {$A$};
\node [below, black] at (61*\xscale,2.9) {$S_i$};
\node [below, black] at (37*\xscale,3.7) {$B$};
\node [right, black] at (112*\xscale,1.5) {$A'$};
\node [right, black] at (100*\xscale,1.5) {$S_j$};
\node [right, black] at (112*\xscale,2.5) {$B'$};
\node [black] at (50*\xscale, 0.8) {$S$};
\end{pgfonlayer}

\draw[black, line width=0.5pt] (42*\xscale, 4) -- (42*\xscale, 6) -- (111*\xscale, 6) -- (111*\xscale, 3);

\draw [black] plot [smooth, tension=1] coordinates { (42*\xscale,2) (42*\xscale,1) (63*\xscale,1.5)  };
\draw [black] plot [smooth, tension=1] coordinates { (100*\xscale,0.5) (110*\xscale,0) (111*\xscale, 1)};

\draw [black, dashed] plot [smooth, tension=1.5] coordinates {(66*\xscale,2) (80*\xscale,1) (98*\xscale,1.5)};

\foreach \i/\j in {2/10, 3/5, 2/27, 1/31, 0/100}{
\pgfmathsetmacro{\maxx}{min(\xscale*\j*2^\i+\xscale*2^\i,\columns*\xscale)}
\draw[neighbor] (\xscale*\j*2^\i,\i) rectangle (\maxx,\i + 1);
}

\foreach \i/\j in {2/10, 1/55, 2/16, 2/17, 3/9, 3/10, 3/11, 2/24, 1/49}{
\pgfmathsetmacro{\maxx}{min(\xscale*\j*2^\i+\xscale*2^\i,\columns*\xscale)}
\draw[active] (\xscale*\j*2^\i,\i) rectangle (\maxx,\i + 1);
}

\foreach \i/\j in {2/10, 1/55, 2/16, 1/49}{
\pgfmathsetmacro{\maxx}{min(\xscale*\j*2^\i+\xscale*2^\i,\columns*\xscale)}
\draw[F] (\xscale*\j*2^\i,\i) rectangle (\maxx,\i + 1);
}

\end{tikzpicture}}
\caption{Proof of Lemma \ref{lem:distance-bound-1}. Left: Step 1. The walk $S$ satisfies (i) and connects $A$ with $A'$. 
If from an inactive box $X$ there is no inactive path to $B$ (dotted line), then there is an active walk $Q$ (dashed line) that 
connects active boxes $E$ and $E'$ in $S$ on either side of $X$. Right: Step 2. The walk $S$ satisfies (i) and contains no inactive 
boxes outside $F$. The boxes $A$, $S_i$, $S_j$ and $A'$ (striped) all belong to $F'$. 
The proof works by finding a path in $F'$ from $S_i$ to $S_j$ (dashed line). 
In both figures active boxes are colored gray and inactive boxes are colored white. \label{fig:lemma8}}
\end{figure}
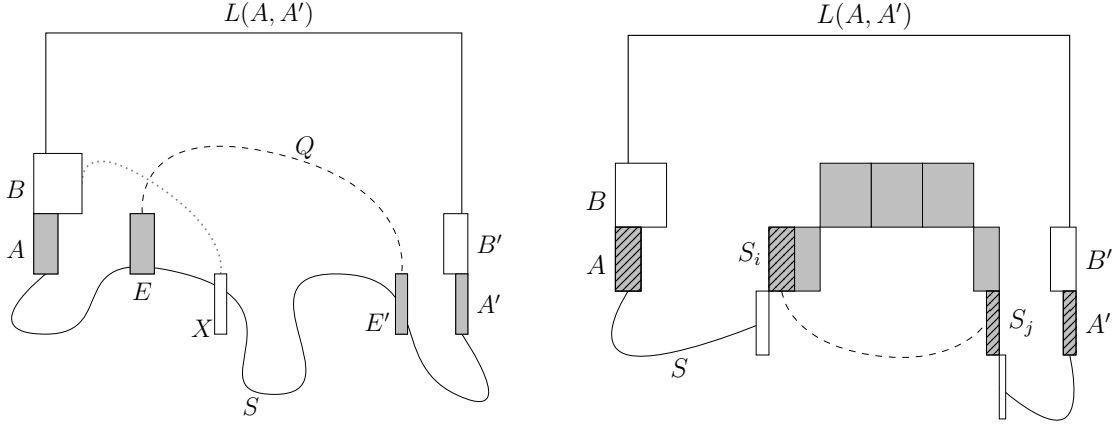

\textit{Step 2.} There is a walk $S$ satisfying (i) that contains no active boxes outside $F'$, where $F'$ is the set of active boxes adjacent 
to a box in $F$. \vspace{0.6\baselineskip}

We start with the walk constructed in Step 1. Since $A$ is adjacent to $B$ it belongs to $F'$. Let $B'$ be the box in $L(A,A')$ directly 
preceding $A'$. We claim that $B'$ belongs to $F$. Note that $B'$ is inactive. 
We use Lemma \ref{obs:barrier} to show that an inactive path from $B$ to $B'$ exists. If such a path would not exist, then an active walk $Q$ would 
exist that intersects all walks from $B$ to $B'$. In particular, $Q$ would contain an active box in $L(A,A') \setminus \{A,A'\}$ 
(which does not lie in $\mathcal{B}_x$, because by assumption $A$ and $A'$ are the only boxes in $L(A,A')$ that belong to $\mathcal{B}_x$) and an 
active box in $S$ (which lies in $\mathcal{B}_x$, because we know there is a path in $\Gamma$ from a vertex in this box to a vertex in $A$). 
This is a contradiction, because by Lemma \ref{lem:B-properties}\eqref{lem:item:B-connections} there cannot be an active walk between a 
box in $\mathcal{B}_x$ and an active box not in $\mathcal{B}_x$. It follows that an inactive path from $B$ to $B'$ exists, so $B'$ belongs 
to $F$. 
Furthermore, every box in $S$ that has an inactive neighbor in $S$ also lies in $F'$, because this inactive neighbor lies in $F$ by Step 1.

Now consider active boxes $S_i$, $S_{i+1}$, \ldots, $S_{j}$ in $S$ such that $S_i$ and $S_j$ lie in $F'$ but 
$S_{i+1}$, \ldots, $S_{j-1}$ do not (Figure \ref{fig:lemma8}, right). 
We claim that there is a path in $F'$ from $S_i$ to $S_j$. Color all boxes in $F'$ blue and all other boxes red. Then our claim is that 
$\mathcal{B}$ contains a blue path from $S_i$ to $S_j$. We use Lemma \ref{obs:barrier} and argue by contradiction. 
If this blue path would not exist, then there would exist a red walk $Q$ that intersects every walk from $S_i$ to $S_j$. 
Because $S_i$ and $S_j$ lie in $F'$, there exists such a walk that apart from $S_i$ and $S_j$ contains only boxes in $F$. 
Because $Q$ does not contain $S_i$ and $S_j$ (which are blue) it must contain a box in $F$. Furthermore, $Q$ also contains one of the active 
boxes $S_{i+1}$, \ldots, $S_j$. Therefore, $Q$ is a connected set of boxes that contains a box in $F$ and an active box. 
This implies that $Q$ must also contain a box in $F'$, which contradicts the fact that $Q$ consists of red boxes. 
This contradiction shows that there must be a blue path in $\mathcal{B}$ from $S_i$ to $S_j$, i.e. a path in $F'$ from $S_i$ to $S_j$. 
We replace the boxes $S_{i+1}$, \ldots, $S_{j-1}$ of $S$ by this path, thereby removing the boxes $S_{i+1}$, \ldots, $S_{j-1}$ from $S$. 
Repeatedly applying this operation, we remove all active boxes that do not lie in $F'$ from $S$. This completes Step 2. \vspace{0.6\baselineskip}

The walk constructed in Step 2 satisfies (i) and (ii'), so we are now done with the case that $L(A,A')$ contains no boxes in 
$\mathcal{B}_x$ other than $A$ and $A'$.

Now suppose $A$ and $A'$ are not the only boxes in $L(A,A')$ that belong to $\mathcal{B}_x$; let $A=A_0$, $A_1$, \ldots, $A_n = A'$ be all 
the boxes in $L(A,A')$ that belong to $\mathcal{B}_x$ (ordered by their position in $L(A,A')$). All these boxes contain vertices in the 
same component of $\Gamma$. For all $i$ we have $L(A_i,A_{i+1}) \subset L(A,A')$ and furthermore $A_i$ and $A_{i+1}$ are the only  boxes 
in $L(A_{i},A_{i+1})$ that belong to $\mathcal{B}_x$. Therefore, a walk from $A_i$ to $A_{i+1}$ satisfying (i) and (ii) exists. 
By concatenating these walks for all $i$ we find a walk $S$ from $A$ to $A'$ satisfying (i) and (ii).

We now construct a path in $\Gamma$ from $x$ to $x'$ of length at most $37|W(A,A')|$. 
We may assume that the active boxes in $S$ are all distinct, because if $S$ contains an active box twice we can remove the intermediate part 
of $S$. 
The number of active boxes in $S$ is at most $9|W(A,A')|$ because each active box in $S$ lies in $W(A,A')$ or is one of the at most $8$ neighbors 
of an inactive box in $W(A,A')$. Suppose $S_i$ and $S_j$ are active boxes in $S$ such that $S_{i+1}$, \ldots, $S_{j-1}$ are all inactive. 
Then for every vertex $v \in S_i$ there is a path in $\Gamma$ of length at most $4$ to a vertex in $S_j$: by (i) there are vertices 
$a \in S_i$, $b \in S_j$ such that $d_{\Gamma}(a,b) \leq 3$ and furthermore $v$ and $a$ are neighbors because they lie in the same box. 
It follows that there is a path of length at most $36|W(A,A')|$ from $x$ to a vertex in $A'$, hence a path of length at 
most $36|W(A,A')|+1 \leq 37|W(A,A')|$ from $x$ to $x'$. This shows that we may take $c=37$.
\end{proof}

\subsection{Bounding the diameter\label{sec:bounding}}

In this subsection we continue with the general setting where $V \subseteq \Ecal_R$ is an arbitrary finite set, and $\Gamma$
is the graph with vertex set $V$ and connection rule $(x,y) \sim (x',y') \iff |x-x'|_{\pi e^{R/2}} \leq e^{\frac12(y+y')}$. 
Here we will translate the bounds from the previous section into results on the maximum diameter of a component of $\Gamma$.
We start with a general observation on graph diameters.

\begin{lem}\label{lem:remove-clique}
Suppose $H_1, H_2$ are induced subgraphs of $G$ such that $V(G) = V(H_1) \cup V(H_2)$ (but $H_1, H_2$ need not be vertex disjoint). 
If every component of $H_1$ has diameter at most $d_1$ and  $H_2$ (is connected and) has diameter at most $d_2$, then every component of $G$ has 
diameter at most $2d_1+d_2+2$.

In particular, if $H_2$ is a clique, then every component of $G$ has diameter at most $2d_1+3$.
\end{lem}

\begin{proof}
Let $C$ be a component of $G$. If $C$ contains no vertices of $H_2$, then $C$ is a component of $H_1$ as well. 
So in this case $C$ has diameter at most $d_1$. 
If $C$ is not a component of $H_1$, then for any vertex $v \in C$ there is a path of length at most $d_1+1$ from $v$ to a vertex in $H_2$. 
Thus, since there is a path of length at most $d_2$ between any two vertices in $H_2$, any two vertices $u,v \in C$ have distance at 
most $(d_1+1) + d_2 + (d_1+1) = 2d_1+d_2+2$ in $G$, as required. 
\end{proof}

We let $\elltil := \lfloor R/2\log 2\rfloor-1$ be the largest $i$ such that layer $i$ is completely below the horizontal line $y = R/2$; 
we set $\Vtil := V \cap \{ y \leq \elltil \cdot \log 2 \}$, and we let $\Gammatil := \Gamma[\Vtil]$ be the subgraph of $\Gamma$ induced by $\Vtil$.
For $A, A' \in \Bcal$ we let $\Wtil = \Wtil(A,A')$ denote the set $W(A,A')$ but corresponding to $\Vtil$ instead of $V$.
(I.e.~boxes in layers $\elltil+1,\dots,\ell$ are automatically inactive. Note that this could potentially increase the size
of $W$ substantially.)

The following lemma gives sufficient conditions for an upper bound on the diameter of each component of $\tilde{\Gamma}$. 
The lemma also deals with graphs that can be obtained by $\tilde{\Gamma}$ by adding a specific type of edges.

\begin{lem}\label{lem:diam-properties}
There exists a constant $c$ such that the following holds.
Let $\Gammatil, \Wtil$ be as above, and let $K = \{(x,y) \in \mathcal{E}_R : y > R/4\}$. 
Consider the following two conditions:
\begin{itemize}
\item[{\bf(i)}] For any two boxes $A$ and $A'$ we have $|\Wtil(A,A')| \leq D$ for some $D$ (possibly depending on $n$);
\item[{\bf(ii)}] There is no inactive path (wrt.~$\Vtil$) in $\mathcal{B}$ connecting a box in $L_0$ with a box in $K$.
\end{itemize}
If {\bf(i)} holds, then each component of $\Gammatil$ has diameter at most $c D$. 
If furthermore {\bf(ii)} holds then, for any any graph $\Gamma'$ that is obtained from $\Gammatil$ by adding an arbitrary set of edges 
$E'$ each of which has an endpoint in $K$, every component of $\Gamma'$ has diameter $c D$.
\end{lem}

\begin{proof}
The first statement directly follows from Lemma~\ref{lem:distance-bound-1}.

If furthermore {\bf(ii)} holds, there exists a cycle of active boxes in $\mathcal{E}_R \backslash K$ that separates $K$ from $L_0$. 
Since vertices in neighboring boxes are connected in $\Gammatil$, this means that there is a cycle in $\Gammatil$ that separates $K$ from $L_0$. 
Every vertex in $K$ lies above some edge in this cycle and thereby lies in the component $C$ of this cycle by Lemma \ref{lem:geometry}(i). 
Thus, every edge of $\Gamma'$ that is not present in $\Gammatil$ has an endpoint in the component $C$ of $\Gammatil$.

Let $d$ be the maximum diameter over all components of $\Gammatil$. An application of Lemma \ref{lem:remove-clique} (with $C$ as one of the two subgraphs; note that we may assume that no added edge connects vertices in the same component, because this can only lower the diameter), we see that the diameter of $\Gamma'$ is at most $3d+2$. This proves the second statement (with a  larger value of $c$).
\end{proof}

\section{Probabilistic bounds}\label{sec:probbounds}

We are now ready to use the results from the previous sections to obtain (probabilistic) bounds on the diameters of components in the KPKVB model. 
Recall from Section \ref{sec:idealized} that $\Gamma_{\alpha,\lambda}$ is a graph with vertex set $V_{\alpha,\lambda}$, where two vertices 
$(x,y)$ and $(x',y')$ are connected by an edge if and only if $|x-x'|_{\pi e^{R/2}} \leq e^{\frac12(y+y')}$. 
Here $V_{\alpha,\lambda}$ is the point set of the Poisson process with intensity 
$f_{\alpha,\lambda} = \mathbf{1}_{\mathcal{E}_R} \lambda e^{-\alpha y}$
on $\mathcal{E}_R = (-\frac{\pi}{2} e^{R/2}, \frac{\pi}{2} e^{R/2}] \times [0,R] \subset \mathbb{R}^2$. 

Consistently with the previous sections,  we define the subgraph $\Gammatil_{\alpha,\lambda}$ of $\Gamma_{\alpha,\lambda}$, induced by the 
vertices in 

\[ \Vtil_{\alpha,\lambda} := \{(x,y) \in V_{\alpha,\lambda} : y \leq (\elltil+1)\log 2\}. \]

\noindent
In the remainder of this section all mention of active and inactive (boxes) will be wrt.~$\Vtil_{\alpha,\lambda}$.

Our plan for the proof of Theorem~\ref{thm:main} is to first show that for $\lambda=\nu\alpha/\pi$ the graph
$\Gammatil_{\alpha,\lambda}$ satisfies the conditions in Lemma~\ref{lem:diam-properties} for some $D = O(R)$. 
In the final part of this section we spell out how this result implies that a.a.s.~all components of the KPKVB random graph $G(N,\alpha,\nu)$ have 
diameter $O(R)$.

We start by showing that property (i) of Lemma \ref{lem:diam-properties} is a.a.s.~satisfied by $\Gammatil_{\alpha,\lambda}$.
To do so, we need to estimate the probability that a box is active if $V$ is the point set $V_{\alpha,\lambda}$ of the Poisson process above.
For $0 \leq i \leq \elltil$, let us write 

\begin{equation}\label{eq:pdef} 
\begin{array}{l} 
p_i = p_{i,\alpha,\lambda} := \Pee_{\alpha,\lambda}( \text{$B$ is active} ), 
\end{array}
\end{equation}

\noindent
where $B \in L_i$ is an arbitrary box in layer $L_i$.

\begin{lem}\label{lem:active-probability}
For each $0 \leq i < \elltil$ we have:

\[ \begin{array}{rcl}  p_i & = & 
1 - \exp\left[ - b \cdot \frac{2^{1-\alpha}}{\alpha} \cdot \lambda \cdot 2^{(1-\alpha)i} \right]  \\
& \geq & 1 - \exp\left[ - \frac{1}{12} \lambda 2^{(1-\alpha)i} \right].
   \end{array}
\]
\end{lem}

\begin{proof}
The expected number of points of $\Vtil_{\alpha, \lambda}$ that fall inside a box $B$ in layer $L_i$ satisfies
\begin{align*}
\mathbb{E}(|B \cap V_{\alpha, \lambda}|) = \int_{i \log(2)}^{(i+1) \log(2)} \int_{0}^{2^i b} \lambda e^{-\alpha y} \dx \dy
= \lambda \cdot b \cdot \frac{1-2^{-\alpha}}{\alpha} \cdot 2^{(1-\alpha)i}
\end{align*}
Since the number of points that fall in $B$ follows a Poisson distribution and 
because $b \cdot \frac{1-2^{-\alpha}}{\alpha} \geq \frac1{12}$, the result follows. 
\end{proof}


\begin{lem}\label{lem:conn}
There exists a $\lambda_0$ such that if $\alpha=1$ and $\lambda > \lambda_0$ then the following holds.
Let $E$ denote the event that there exists a connected subgraph $C \subseteq \Bcal$ with $|C| > R$ such that least half of the boxes of $C$ are inactive.
Then $\Pee_{1,\lambda}(E) = O(N^{-1000})$.
\end{lem}

\begin{proof}
The proof is a straightforward counting argument.
If $C$ is a connected subset of the boxes graph $\Bcal$ and $A \in C$ is a box of $C$, then
there exists a walk $P$, starting at $A$, through all boxes in $C$, that uses no edge in $\mathcal{B}$ more than twice 
(this is a general property of a connected graph). Since the maximum degree of $\Bcal$ is 8, the walk $P$ visits no box more than 8 times.
Thus, the number of connected subgraphs of $\Bcal$ of cardinality $i$ is no more than $|\Bcal| \cdot 8^{8i} = (2^{\ell+1}-1) \cdot 2^{24 i}
= e^{O(R)} \cdot 2^{24 i}$ (using the definition~\eqref{eq:elldef} of $\ell$).
Given such a connected subgraph $C$ of cardinality $i>R$ there are $\binom{i}{i/2}$ ways to choose a subset of cardinality $i/2$. 
Out of any such subset at most 63 boxes lie above level $\elltil$ by Lemma~\ref{lem:B-properties} and by 
Lemma~\ref{lem:active-probability} each of the remaining $i/2-63 > i/4$ is inactive with probability at most $e^{-\lambda/12}$.
This gives

\[ \begin{array}{rcl} 
\Pee_{1,\lambda}( E ) 
&  \leq &  \sum_{i>R} |\Bcal| \cdot 8^{8i} \cdot \binom{i}{i/2} \cdot e^{-(i/2-63) \cdot \lambda / 12} \\
& \leq &  e^{O(R)} \cdot \sum_{i>R} 8^{8i} \cdot 2^i \cdot e^{-{i/4} \cdot \lambda / 12} \\
& = & e^{O(R)} \cdot O\left( \left( 2^{25} \cdot e^{-\lambda/48} \right)^R \right) \\
& = & \exp[ O(R) - \lambda \cdot \Omega(R) ] \\
& = & O( N^{-1000} ),
\end{array}
\]

\noindent
where the third and fifth line follow provided $\lambda$ is chosen sufficiently large.
\end{proof}

\begin{cor}\label{cor:diameter-gamma-1}
There exist constants $c, \lambda_0$ such that if $\alpha = 1$ and $\lambda>\lambda_0$
then

\[ \Pee_{1,\lambda}\left( \text{there exist boxes $A, A'$ with $|\Wtil(A,A')|> c R$} \right) = O( N^{-1000} ). \]

\end{cor}

\begin{proof}
We let $\lambda_0$ be as provided by Lemma~\ref{lem:conn} and we take $c := 5$.
We note that for every two boxes $A, A'$ the set $\Wtil(A,A')$ is a connected set, and all boxes
except for some of the at most $2R$ boxes on $L(A,A')$ must be inactive by definition of 
$\Wtil$.
Hence if it happens that $|\Wtil(A,A')| > 5R$ for some pair of boxes $A, A'$ then event $E$ defined in Lemma~\ref{lem:conn} holds.
The corollary thus follows directly from Lemma~\ref{lem:conn}. 
\end{proof}

We now want to show that in the case when $\frac12 < \alpha < 1$ and $\lambda > 0$ we also have that, with probability very close to one, 
$|\Wtil(A,A')| = O(R)$ holds for all $A, A'$.
Recall that the probability that a box in layer $i$ is inactive is upper bounded by $\exp(-\frac{\lambda}{12} 2^{(1-\alpha)i})$ 
(this bound now depends on $i$), which decreases rapidly if $i$ increases. 
However, for small values of $i$ this expression could be very close to $1$, depending on the value of $\lambda$. 
In particular we cannot hope for something like Lemma~\ref{lem:conn} to hold for all $\frac12 < \alpha < 1$ and $\lambda > 0$.

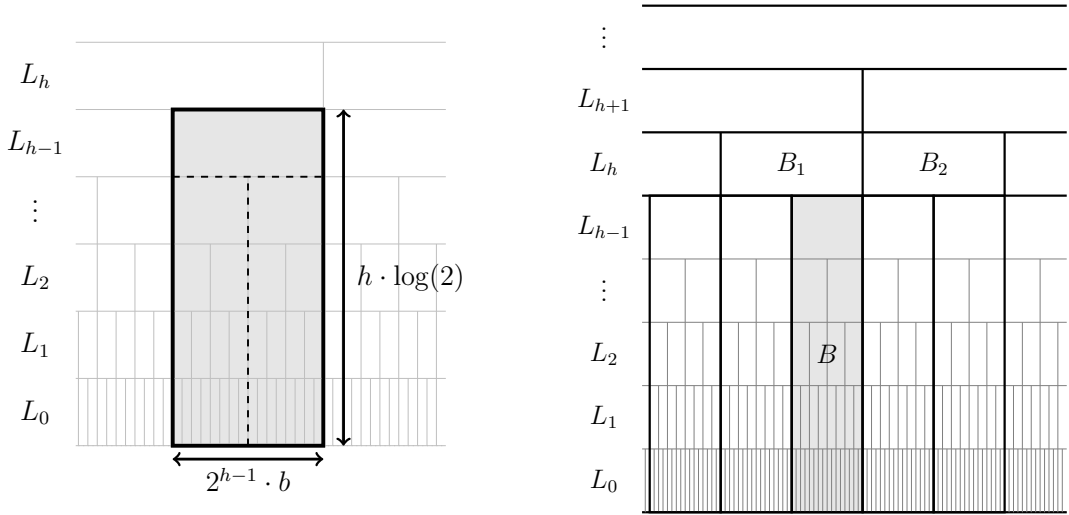
\begin{figure}[h!]
\centering
\resizebox{!}{0.38\linewidth}{
\begin{tikzpicture}

\pgfdeclarelayer{bg}    
\pgfdeclarelayer{edges}
\pgfdeclarelayer{text}
\pgfsetlayers{bg,edges,main,text}

\def\xscale{0.07}
\def\columns{90}
\def\rows{8}
\def\density{0.5}
\def\rowsh{4}
\pgfmathtruncatemacro{\rowsminusone}{\rows - 1}
\pgfmathtruncatemacro{\rowsplusone}{\rows + 1}
\pgfmathtruncatemacro{\rowsplustwo}{\rows + 2}
\pgfmathtruncatemacro{\columnsa}{\columns/2}

\def\yshift{0}

\tikzset{%
	examplebox/.style = {color=gray!70!white},
    neighbor/.style={color=red,line width=1pt, fill=red!10!white},
}
\draw [white] (0,0) -- (1,0);

\draw [black, line width=1.5pt,fill=black!10!white]  (32*\xscale,\yshift+1) rectangle (64*\xscale,\yshift+6);


\foreach \x in {1,...,\rowsminusone}{
\draw [gray!50!white] (0.8,\yshift+ \x) -- (\columns*\xscale, \yshift+\x);
}


\foreach \y in {0,...,\rowsminusone} {
    \pgfmathtruncatemacro{\columnsa}{\columns/2^\y}
	\foreach \x in {0,...,\columnsa}{
	\pgfmathsetmacro{\z}{2^\y*\x} 
	\pgfmathtruncatemacro{\w}{\z}
	\pgfmathtruncatemacro{\u}{\z*\xscale*5}
	\ifthenelse{\w < \columns}{
		\ifthenelse{\u>3}{
		\pgfmathsetmacro{\v}{max(\y,1)};
		\draw [gray!50!white] (\z*\xscale,\yshift+\v) -- (\z*\xscale,\yshift+\y+1);
		}{}
	}{}
	} 
}

\draw [black, line width=1.5pt]  (32*\xscale,\yshift+1) rectangle (64*\xscale,\yshift+6);

\draw [black,<->,line width=1.2pt] (32*\xscale,\yshift+0.8) -- (64*\xscale,\yshift+0.8);
\draw [black,<->,line width=1.2pt] (64*\xscale + 0.3, \yshift+1) -- (64*\xscale + 0.3,\yshift+6);

\node [black] at (48*\xscale,\yshift+0.5)  {$2^{h-1} \cdot b$};
\node [black] at (64*\xscale + 1.3,\yshift+3.5) {$h \cdot \log(2)$};

\draw [black,line width=0.8pt,dashed] (48*\xscale,1) -- (48*\xscale,5);
\draw [black,line width=0.8pt,dashed] (32*\xscale,5) -- (64*\xscale,5);

\foreach \y in {0,1,2}{
\node [black] at (0.2,\y+1.5) {$L_\y$};
}

\node [black] at (0.2,4.6) {$\vdots$};
\node [black] at (0.2,5.5) {$L_{h-1}$};
\node [black] at (0.2,6.5) {$L_h$};

	\end{tikzpicture}}\hspace{1cm}
\resizebox{!}{0.41\linewidth}{
\begin{tikzpicture}

\pgfdeclarelayer{bg}    
\pgfdeclarelayer{edges}
\pgfdeclarelayer{text}
\pgfsetlayers{bg,edges,main,text}

\def\xscale{0.07}
\def\columns{110}
\def\rows{8}
\def\density{0.5}
\def\rowsh{4}
\pgfmathtruncatemacro{\rowsminusone}{\rows - 1}
\pgfmathtruncatemacro{\rowsplusone}{\rows + 1}
\pgfmathtruncatemacro{\rowsplustwo}{\rows + 2}
\pgfmathtruncatemacro{\columnsa}{\columns/2}

\tikzset{%
	examplebox/.style = {color=gray!70!white},
    neighbor/.style={color=black,line width=1pt},
}

\draw [fill= gray!20!white, line width=0pt] (48*\xscale+0.01,0.01) rectangle (64*\xscale-0.01,5-0.01);



\foreach \x in {0,...,\rows}{
\ifthenelse{\x > 4}{
\draw [black, line width=1pt] (1, \x) -- (\columns*\xscale, \x);
}{
\draw [gray] (1, \x) -- (\columns*\xscale, \x);}
}

\draw [black, line width=1pt] (1,0) -- (\columns*\xscale,0);


\foreach \y in {0,...,\rowsminusone} {
    \pgfmathtruncatemacro{\columnsa}{\columns/2^\y}
	\foreach \x in {0,...,\columnsa}{
	\pgfmathsetmacro{\z}{2^\y*\x} 
	\pgfmathtruncatemacro{\w}{\z}
	\pgfmathtruncatemacro{\v}{\z*\xscale}
	\ifthenelse{\w < \columns}{
	\ifthenelse{\y > 4}{
	\ifthenelse{\v>0}{
	\draw [black, line width=1pt] (\z*\xscale,\y) -- (\z*\xscale,\y+1);}{}
	}{
	\ifthenelse{\v>0}{
	\draw [gray] (\z*\xscale,\y) -- (\z*\xscale,\y+1);}}{}
	}{}
	} 
}

\foreach \i/\j in {0/16,0/32,0/48,0/64,0/80}{
\begin{pgfonlayer}{main}
\pgfmathsetmacro{\maxx}{min(\xscale*\j*2^\i+16*\xscale*2^\i,\columns*\xscale)}
\draw[neighbor] (\xscale*\j*2^\i,\i) rectangle (\maxx,\i + 5);
\end{pgfonlayer}
}


\foreach \y in {0,1,2}{
\node [black] at (1-0.6,\y+0.5) {$L_\y$};
}

\node [black] at (1-0.6,3.6) {$\vdots$};
\node [black] at (1-0.6,4.5) {$L_{h-1}$};
\node [black] at (1-0.6,5.5) {$L_h$};
\node [black] at (1-0.6,6.5) {$L_{h+1}$};
\node [black] at (1-0.6,7.6) {$\vdots$};


\node [black] at (56*\xscale,2.5) {$B$};
\node [black] at (48*\xscale,5.5) {$B_1$};
\node [black] at (80*\xscale,5.5) {$B_2$};

\end{tikzpicture}}
\caption{Left: an $h$-block (in the figure $h=5$). An $h$-block is the union of $2^h-1$ boxes in the lowest $h$ layers. 
Right: definition of a lonely block (used in the proof of Lemma \ref{lem:diameter-gamma}). The lowest $h$ layers are partitioned 
into $h$-blocks. 
An $h$-block $B$ in $W'$ is called \emph{lonely} if both boxes $B_1$ and $B_2$ lying above it are not in $W'$. If $|W'| > 3$ 
and $B$ is lonely, one of the blocks adjacent to $B$ contains a horizontal path in $W$. \label{fig:hblock}}
\end{figure}

To gain control over the boxes in the lowest layers, we merge boxes in the lowest layers into larger blocks. 
An \emph{$h$-block}\index{block}\index{$h$-block} is defined as the union of a box in $L_{h-1}$ and all $2^h-2$ boxes lying below 
this box (Figure \ref{fig:hblock}, left). In other words, an $h$-block consists of $2^h-1$ boxes in the lowest $h$ layers that 
together form a rectangle. The following Lemma shows that the probability that an $h$-block contains a horizontal inactive path 
can be made arbitrarily small by taking $h$ large.

Let us denote by $q_h$ the probability:

\begin{equation}\label{eq:cdef} 
q_h = q_{h,\alpha,\lambda} := \Pee_{\alpha,\lambda}\left(  \text{$H$ has a vertical, active crossing}\right),
\end{equation}

\noindent
where $H$ is an arbitrary $h$-block; and a ``vertical, active crossing'' means a path of active boxes (in $\Bcal^\ast$) inside the block connecting
the unique box in the highest layer to a box in the bottom layer.

\begin{lem}\label{lem:ph}\index{$h$-block}\index{block}\index{box}
If $\alpha < 1$ and $\lambda > 0$ then, for every $\eps>0$, there exists an $h_0=h_0(\eps,\alpha,\lambda)$ such that $q_h > 1-\eps$ for 
all $h_0 \leq h \leq \elltil$.
\end{lem}

\begin{proof}
In the proof that follows, we shall always consider blocks that do not extend above the horizontal line $y=R/2$, 
(i.e.~boxes in layers $h \leq R/2\log 2-1$) so that we can use Lemma \ref{lem:active-probability} to estimate the 
probability that a box is active. 

An $(h+1)$-block $H$ consists of one box $B$ in layer $L_h$ and two $h$-blocks $H_1$, $H_2$.
There is certainly a vertical, active crossing in $H$ if $B$ is active and either $H_1$ or $H_2$ 
has a vertical, active crossing.
In other words,

\begin{equation}\label{eq:aap} 
q_{h+1} \geq p_h \cdot (2 q_h - q_h^2), 
\end{equation}

\noindent
where $p_h \geq 1 - \exp[ - \frac{1}{12} \lambda 2^{(1-\alpha)h} ]$ is the probability that $B$ is active.
We choose $\delta = \delta(\eps)$ small, to be made precise shortly.
Clearly there is an $h_0$ such that $p_h > 1-\delta$ for all $h_0 \leq h \leq \elltil$.
Thus~\eqref{eq:aap} gives that $q_{h+1} \geq f(q_h)$ for all such $h$, where $f(x) = (1-\delta)(2x-x^2)$.
It is easily seen that $f$ has fixed points $x = 0, \frac{1-2\delta}{1-\delta}$,
that $x < f(x) < \frac{1-2\delta}{1-\delta}$ for $0 < x < \frac{1-2\delta}{1-\delta}$ and $\frac{1-2\delta}{1-\delta} < f(x) < x$ 
for $\frac{1-2\delta}{1-\delta} < x \leq 1$.
Therefore, using that clearly $0 < q_{h_0} < 1$ (there is for instance a strictly positive probability all boxes of the block $H$ are 
active, resp.~inactive),
we must have $f^{(k)}(q_{h_0} ) \to \frac{1-2\delta}{1-\delta}$ as $k\to\infty$, where $f^{(k)}$ denotes the $k$-fold composition of $f$ 
with itself.
Hence, provided we chose $\delta=\delta(\eps)$ sufficiently small, there is a $k_0=k_0(\eps)$ such that
$q_{h_0+k} \geq f^{(k)}(q_{h_0}) > 1-\eps$ for all $k_0 \leq k \leq \elltil - h_0$. 
\end{proof}

\begin{lem}\label{lem:diameter-gamma}
For every $\alpha < 1$ and $\lambda > 0$ there exists a $c = c(\alpha,\lambda)$ such that

\[ \Pee_{\alpha,\lambda}\left( \text{there exist boxes $A$ and $A'$ with $|\Wtil(A,A')| > c R$} \right) = O( N^{-1000} ).  
\]

\end{lem}

\begin{proof}
Let $p_i$ be as defined in~\eqref{eq:pdef} and $q_i$ as defined in~\eqref{eq:cdef}.
Let $\eps>0$ be arbitrary, but fixed, to be determined later on in the proof.
By Lemmas~\ref{lem:active-probability} and~\ref{lem:ph}, there exists an $h$ such that 

\[ p := \min\{ q_h, p_i : h \leq i \leq \elltil \} > 1-\eps. \]

\noindent
We now create a graph $\Bcal'$, modified from  the boxes graph $\Bcal$, as follows.
The vertices of $\Bcal'$ are the boxes above layer $h$, together with the $h$-blocks.
Boxes or blocks are neighbors in $\Bcal'$ if they share at least a corner. 
Note that the maximum degree of $\mathcal{B}'$ is at most $8$. 

Given two boxes $A, A' \in \Bcal$ we define $W'(A,A') \subseteq \Bcal'$ as the natural analogue of $\Wtil(A,A')$, 
i.e.~the set of all boxes of $\Wtil(A,A')$ above layer $h$ together with all $h$-blocks that contain at least one element of $\Wtil(A,A')$.
Note that $W'(A,A')$ is a connected set in $\Bcal'$ and that $|W'(A,A')| \geq |\Wtil(A,A)| / (2^h-1)$.

We will say that an $h$-block $B$ is \emph{lonely} if the 
two boxes in $L_h$ adjacent to $B$ both do not lie in $\Wtil(A,A')$ (Figure \ref{fig:hblock}, right).
Observe that if $B$ is lonely, then at least one of the two neighbouring blocks must have a horizontal, inactive crossing.
This shows that:

\begin{equation}\label{eq:obss} |\text{blocks without an active, vertical crossing}| \geq |\text{lonely blocks}|/2. \end{equation}

\noindent
Consider two boxes $A, A' \in \Bcal$ and assume that $|\Wtil(A,A')| > c R$, where $c$ is a large constant to be made precise later.
By a previous observation $|W'(A,A')| \geq (\frac{c}{2^h-1}) R =: d R$.
We distinguish two cases. 

{\bf Case a):} at least $|W'(A,A')| / 100$ of the elements of $W'(A,A')$ are boxes (necessarily above layer $h$).
Subtracting the at most 63 boxes of levels $\elltil+1,\dots, \ell$ and the at most $2R$ boxes of $L(A,A')$, we see that 
at least $|W'(A,A)| / 100 - (2R + 63) \geq |W'(A,A')| / 1000$ boxes of $W'(A,A')$ must be inactive and lie in levels $h, \dots, \elltil$.
(Here the inequality holds assuming $|W'(A,A')| \geq d R$ with $d$ sufficiently large.) 

{\bf Case b):} at most $|W'(A,A')| / 100$ of the elements of $W'(A,A')$ are boxes. Hence, at least 
$\frac{99}{100} |W'(A,A')|$ of the elements of $W'(A,A')$ are $h$-blocks.
Of these, at least $\frac{97}{100} |W'(A,A')|$ blocks must be lonely, since each box of $W'$ is adjacent to no more than 
two $h$-blocks of $W'$. Thus, by the previous observation~\eqref{eq:obss}
at least $\frac{97}{200} |W'(A,A')| \geq |W'(A,A')| / 1000$ elements of $W'$ are blocks without a vertical, active crossing.

Combining the two cases, we see that either $W'$ contains $|W'(A,A')|/1000$ inactive boxes in the levels $h$, \ldots, $\elltil$, or $W'$ contains  $|W(A,A')|/1000$ blocks without a vertical, active crossing. Summing over all possible choices of $A,A'$ and all possible sizes of $W'(A,A')$, we see that

\[ 
 \begin{array}{rcl} 
\Pee_{\alpha,\lambda}\left( \text{there exist $A, A'$ with $|\Wtil(A,A')|> c R$} \right)
& \leq & |\Bcal|^2  \cdot \sum_{i\geq d R} 8^{8i} \cdot (1-p)^{i/1000} \\ 
& \leq &  |\Bcal|^2 \cdot \sum_{i\geq d R} \left( 8^8 \cdot \eps^{1/1000} \right)^i \\
& = &  |\Bcal|^2 \cdot O\left( \left( 8^8 \cdot \eps^{1/1000} \right)^{dR} \right) \\
& = & \exp\left[O(R) - d \cdot \Omega(R) \right] \\ 
& = & O(N^{-1000}),
\end{array}
\]

\noindent
where the factor $8^{8i}$ in the first line is a bound on the number of connected subsets 
of $\Bcal'$ of cardinality $i$ that contain $A, A'$; the third line holds provided $\eps$ is sufficiently small ($\eps < 8^{-8000}$ will do); 
and the last line holds provided $c$ (and thus also $d = c/(2^h-1)$) was chosen sufficiently large.
\end{proof}

We now turn to the proof of (ii) of Lemma \ref{lem:diam-properties}. 

\begin{lem}\label{lem:property2}
If either
 \begin{enumerate} [(i)]
\item $\frac12 < \alpha < 1$ and $\lambda > 0$ is arbitrary, or; 
\item $\alpha = 1$ and $\lambda$ is sufficiently large,
 \end{enumerate} \index{KPKVB model!diameter}
then, it holds with probability $1-O(N^{-1000})$ that there are no inactive paths in $\mathcal{B}$ from $L_0$ to 
$K :=  \{(x,y) \in \mathcal{E}_R : y > R/4\}$. 
\end{lem}

\begin{proof}
Since only the boxes below the line $y=R/2$ are relevant, we can freely use Lemma \ref{lem:active-probability}. 
Note that an inactive path in $\mathcal{B}$ from $L_0$ to $K$ would have length at least $R/4$ (the height of each layer equals $\log 2 < 1$) 
and that it would have a subpath of length at least $R/8$ that lies completely in $\{(x,y) : y > R/8\}$. Let $q$ be the maximum probability 
that a box between the lines $y=R/8$ and $y=R/2$ is inactive. Since there are $\exp(O(R))$ boxes and at most $9^k$ paths of length $k$ starting 
at any given box, the probability that such a subpath exists is at most $\exp(O(R)) 9^{R/8} q^{R/8} = \exp(O(R) + \log(q) R/8)$. 
If $\alpha=1$ then $q \leq \exp(-\lambda/12)$, which can be chosen arbitrarily small by choosing $\lambda$ sufficiently large. 
For sufficiently small $q$ we then have $\exp(O(R) + \log(q) R/8) \leq \exp(-R/2) = O(N^{-1000})$ and therefore such a path does 
not exist with probability $1- O(N^{-1000})$. 
If $\alpha<1$ we have $q \leq \exp(-\lambda/12 \cdot 2^{(1-\alpha)R/8})$ and it follows 
that $\exp(O(R) + \log(q) R/8) = \exp(O(R) - \lambda/12 \cdot 2^{(1-\alpha)R/8} \cdot R/8) = \exp(-\omega(R) )$, so we 
can draw the same conclusion. 
\end{proof}

We are almost ready to finally prove Theorem~\ref{thm:main}, but it seems helpful to 
first prove a version of the theorem for $G_{\Po}$, the Poissonized version of the model.

\begin{proposition}[Theorem \ref{thm:main} for $G_{\Po}$]
If either
\begin{enumerate}[(i)]
\item $\frac12 < \alpha < 1$ and $\nu > 0$ is arbitrary, or; 
\item $\alpha = 1$ and $\nu$ is sufficiently large,
\end{enumerate}\index{KPKVB model!diameter}
then, a.a.s., every component of $G_{\Po}(N; \alpha, \nu)$ has diameter $O(\log(N))$.
\end{proposition}

\begin{proof}\index{coupling!between KPKVB model and idealized model}
Let $\Gtil_{\Po}$ be the subgraph of $G_{\Po}$ induced by the vertices of radius larger than $R - (\elltil+1)\log 2$. 
(Here $\elltil := \lfloor R/2\log 2\rfloor - 1$ is as before.) 
By the triangle inequality, all vertices of $G_{\Po}$ with distance at most $R/2$ from the origin form a clique. 
Moveover, by Lemmas~\ref{lem:coupling},~\ref{lem:coupling-edges} and~\ref{lem:B-properties}, a.a.s.~the vertices of $G_{\Po}$ with radii
between $R/2$ and $R-(\elltil+1)\log 2$ can be partitioned into up to 63 cliques corresponding to the boxes
above level $\elltil$. 
In other words, a.a.s., $\Gtil_{\Po}$ can be obtained from $G_{\Po}$ by successively removing up to 64 cliques.
Therefore, by up to 64 applications of Lemma \ref{lem:remove-clique}, it suffices to show that a.a.s.~every component of 
$\Gtil_{\Po}$ has diameter $O(\log N)$. 
Again invoking Lemmas~\ref{lem:coupling} and~\ref{lem:coupling-edges} as well as Lemma~\ref{lem:diam-properties}, it thus suffices 
to show that a.a.s.~$\Gammatil_{\alpha,\nu\alpha/\pi}$ satisfies the conditions (i) and (ii) of Lemma \ref{lem:diam-properties}. 
This is taken care of by Corollary~\ref{cor:diameter-gamma-1} and Lemma~\ref{lem:property2} in the case when $\alpha=1$ and $\nu$ is 
sufficiently large and Lemmas~\ref{lem:diameter-gamma} and~\ref{lem:property2} in the case when $\alpha < 1$ and $\nu>0$ is arbitrary.
\end{proof}

Finally, we are ready to give a proof of Theorem~\ref{thm:main}.

\begin{proofof}{Theorem~\ref{thm:main}}
Let us point out that $G_{\Po}$ conditioned on $Z=N$ has exactly the same distribution as $G = G(N; \alpha,\nu)$.
We can therefore repeat the previous proof, where we substitute the use of Lemmas~\ref{lem:coupling} and~\ref{lem:coupling-edges} 
by Corollary~\ref{cor:koppel}, but with one important additional difference.
Namely, now we have to show that $\Gammatil_{\alpha,\nu\alpha/\pi}$ satisfies the conditions (i) and (ii) of Lemma \ref{lem:diam-properties} 
a.a.s.~{\bf conditional on $Z=N$}.

To this end, let $E$ denote the event that $\Gammatil_{\alpha,\nu \alpha/\pi}$ fails to have one or both of properties (i) or (ii). 
By Corollary~\ref{cor:diameter-gamma-1}, resp.~Lemma~\ref{lem:diameter-gamma}, and Lemma~\ref{lem:property2} we have $\mathbb{P}(E) = O(N^{-1000})$. 
Using the standard fact that $\Pee( \Po(N)=N ) = \Omega( N^{-1/2} )$, it follows that 

\[ \mathbb{P}( E \mid Z=N) \leq \mathbb{P}(E)/\mathbb{P}(Z=N) = O(N^{-1000})/\Omega(N^{-1/2}) = o(1), \]

\noindent
as required.
\end{proofof}

\section{Discussion and further work\label{sec:conclusion}}

In this paper we have given an upper bound of $O(\log N )$ on the diameter of the components of the KPKVB random graph, which holds  
when $1/2 < \alpha < 1$ and $\nu > 0$ is arbitrary and when $\alpha=1$ and $\nu$ is sufficiently large.
Our upper bound is sharp up to the leading constant hidden inside the $O(\cdot)$-notation.

The proof proceeds by considering the convenient idealized model introduced by Fountoulakis and the first author~\cite{FM}, and a 
relatively crude discretization of this idealized model. 
The discretization is obtained by dissecting the upper half-plane into rectangles (``boxes'') and declaring a box active 
if it contains at least one point of the idealized model. With a mix of combinatorial and geometric arguments, we are then 
able to give a deterministic upper bound on the component sizes of the idealized model in terms of the combinatorial 
structure of the active set of boxes, and finally we apply Peierls-type arguments to give a.a.s.~upper bounds for the diameters of
all components of the idealized model and the KPKVB model.

We should also remark that the proof in~\cite{BFM_connected} that $G(N;\alpha,\nu)$ is a.a.s.~connected when $\alpha < 1/2$ in fact shows
that the diameter is a.a.s.~$O(1)$ in this case.
What happens with the diameter when $\alpha=1/2$ is an open question. 

What happens when $\alpha > 1$ or $\alpha=1$ and $\nu$ is arbitrary is another open question.
In the latter case our methods seem to break down at least partially. 
We expect that the relatively crude discretization we used in the current paper will not be helpful here and a more refined proof 
technique will be needed.

Another natural question is to see whether one can say something about the difference between the diameter of the largest component and 
the other components. We remark that by the work of Friedrich and Krohmer~\cite{FK} the largest component
has diameter $\Omega( \log N )$ a.a.s., but their argument can easily be adapted to show that there will be 
components other than the largest component that have diameter $\Omega(\log N )$ as well.

Another natural, possibly quite ambitious, goal for further work would be to determine the leading constant (i.e.~a constant $c=c(\alpha,\nu)$ such
that the diameter of the largest component is $(c+o(1)) \log N$ a.a.s.) if it exists, or indeed even just to establish the existence of a 
leading constant without actually determining it.
We would be especially curious to know if anything special happpens as $\nu$ approaches $\nu_{\text{c}}$.

We have only mentioned questions directly related to the graph diameter here. To the best of our knowledge the study of (most) 
other properties of the KPKVB model is largely virgin territory.

\section*{Acknowledgements}

We warmly thank Nikolaos Fountoulakis for helpful discussions during the early stages of this project.

\bibliographystyle{plain}

\bibliography{References}

\appendix

\section{The proof of Corollary~\ref{cor:koppel}\label{sec:appendix}}

We show that Lemma \ref{lem:coupling} and \ref{lem:coupling-edges} are also true when conditioned on $Z=N$.
Recall that we say that an event $A$ happens a.a.s. conditional on $B$ if $\mathbb{P}(A \mid B) \to 1$ as $N \to \infty$.

\begin{lem}[Lemma \ref{lem:coupling} conditional on $Z=N$] \label{lem:coupling-2}
Let $\alpha>\frac12$. On the coupling space of Lemma \ref{lem:coupling}, conditional on $Z=N$, a.a.s.~$V_{\nu \alpha/\pi}$ is the 
image of the vertex set of $G_{\Po}$ under $\Psi$. 
\end{lem}

\begin{proof}
Write $V = \{X_1,\ldots,X_Z\}$ and $\tilde{V} = V_{\nu \alpha/\pi}$. As in $\cite{FM}$, there are independent Poisson processes 
$\mathcal{P}_0$, $\mathcal{P}_1$, $\mathcal{P}_2$ on $\mathcal{E}_R$ such that $\Psi(V) = \mathcal{P}_0 \cup \mathcal{P}_1$, 
$\tilde{V} = \mathcal{P}_0 \cup \mathcal{P}_2$ and $\mathbb{E}|\mathcal{P}_1|$, $\mathbb{E}|\mathcal{P}_2| = o(1)$. 
We now find

\begin{align*}
\mathbb{P}(\tilde{V} = \Psi(V) \mid Z=N) &= \mathbb{P}(|\mathcal{P}_1| = |\mathcal{P}_2| = 0 \mid |\mathcal{P}_0| + |\mathcal{P}_1| = N) \\
&= \mathbb{P}(|\mathcal{P}_1| =0\mid |\mathcal{P}_0| + |\mathcal{P}_1| = N) \mathbb{P}(|\mathcal{P}_2| = 0)
\end{align*}

because $\mathcal{P}_0$, $\mathcal{P}_1$ and $\mathcal{P}_2$ are independent. From $\mathbb{E}|\mathcal{P}_2|=o(1)$ it follows 
that $\mathbb{P}(|\mathcal{P}_2|=0) = 1-o(1)$. Furthermore, since the conditional distribution of a Poisson distributed variable 
given its sum with an independent Poisson distributed variable is binomial, we have
$
\mathbb{P}(|\mathcal{P}_1| =0\mid |\mathcal{P}_0| + |\mathcal{P}_1| = N) 
= \binom{N}{N} (1 - \mathbb{E}|\mathcal{P}_1|/N )^N = ( 1 - o(1)/N )^N = 1-o(1),
$
from which it follows that $\mathbb{P}(\tilde{V} = \Psi(V) \mid Z=N) = 1-o(1)$.
\end{proof} 

\begin{lem}[Lemma \ref{lem:coupling-edges} conditional on $Z=N$]\label{lem:coupling-edges-2}
Let $\alpha>\frac12$. On the coupling space of Lemma \ref{lem:coupling}, conditional on $Z=N$, a.a.s. it holds for $1 \leq i,j \leq Z$ that
\begin{itemize}
\item[(i)] if $r_i,r_j \geq \frac12R$ and $\tilde{X}_i\tilde{X}_j \in E(\Gamma_{\alpha,\nu \alpha/\pi})$, then $X_iX_j \in E(G_{\Po})$.
\item[(ii)] if $r_i,r_j \geq \frac34 R$, then $\tilde{X}_i\tilde{X}_j \in E(\Gamma_{\alpha,\nu \alpha/\pi}) \iff X_iX_j \in E(G_{\Po})$. 
\end{itemize} 
Here $r_i$ and $r_j$ denote the radial coordinates of $X_i,X_j \in \mathcal{D}_R$.
\end{lem}

\begin{proof}
Let $A$ denote the event that (i) or (ii) fails for some $i,j \leq Z$ and let $B$ denote the event that (i) or (ii) fails for 
some $i,j \leq \min(N,Z)$. It follows that

\begin{equation}\label{eq:app1}
\mathbb{P}(B \mid Z \geq N) \leq \frac{\mathbb{P}(B)}{\mathbb{P}(Z \geq N)} \leq \frac{\mathbb{P}(A)}{\mathbb{P}(Z \geq N)} 
\stackrel{N \to \infty}{\longrightarrow} \frac{0}{1/2} = 0,
\end{equation}

\noindent
because $\mathbb{P}(A) \to 0$ by Lemma \ref{lem:coupling-edges} and $\mathbb{P}(Z \geq N) \to \frac12$ by the central limit theorem. 
Next, let us observe that $\Pee( B \mid Z = N ) = \Pee( B \mid Z=N+1 ) = \dots$ since the points with index greater than 
$\min(N,Z)$ are irrelevant for the event $B$. 
From this it follows that 

\begin{equation}\label{eq:app2}
\Pee( B \mid Z \geq N ) = \frac{\sum_{i \geq N} \Pee( B \mid Z = i )\Pee( Z=i )}{\sum_{i\geq N} \Pee( Z = i ) } = \Pee(B\mid Z=N ).
\end{equation}

\noindent
Combining~\eqref{eq:app1} and~\eqref{eq:app2} we see that $\Pee( B \mid Z = N ) = o(1)$, as desired.
\end{proof}

\end{document}